\documentclass[11pt,reqno]{amsart}
\usepackage{xy,epsfig,amsfonts,amssymb,amsmath,epstopdf}

\def\rto{\rightarrow}

\def\S{{\mathbb S}}

\def\PP{{\mathcal P}}
\def\Sym{{\mathcal S}}

\DeclareMathOperator{\kk}{\mathbf{k}}
\DeclareMathOperator{\un}{\mathbf{1}}

\DeclareMathOperator{\Com}{\mathcal{C}om}
\DeclareMathOperator{\Zin}{\mathcal{Z}in}
\DeclareMathOperator{\Dend}{\mathcal{D}end}
\DeclareMathOperator{\TD}{\mathcal{T}riDend}
\DeclareMathOperator{\CTD}{\mathcal{C}TD}
\DeclareMathOperator{\Comp}{\mathcal{C}omp}
\DeclareMathOperator{\As}{\mathcal{A}s}
\DeclareMathOperator{\Lie}{\mathcal{L}ie}
\DeclareMathOperator{\Tree}{\mathcal{T}ree}

\def\twb{twisted bialgebra}

\DeclareMathOperator{\st}{st}
\DeclareMathOperator{\Id}{Id}
\DeclareMathOperator{\Ker}{Ker}

\DeclareMathOperator{\Prim}{Prim}
\DeclareMathOperator{\Sh}{Sh}

\DeclareMathOperator{\CC}{\mathbf{C}}
\DeclareMathOperator{\Oubli}{\mathcal{O}}

\DeclareMathOperator{\grVect}{\mathbf{grVect}}
\DeclareMathOperator{\Coalg}{\mathbf{Coalg}}
\DeclareMathOperator{\Smod}{\mathbf{S-mod}}
\DeclareMathOperator{\Smodr}{\mathbf{S-mod_r}}

\xyoption{all}
\CompileMatrices

\begin{document}

\title[From left modules to algebras]{From left modules to algebras over an operad: application to combinatorial hopf algebras}
\author[M. Livernet]{Muriel Livernet}
\address{Institut Galil\'ee\\
Universit\'e Paris Nord\\
  93430 Villetaneuse       \\
         France}
\email{livernet@math.univ-paris13.fr}
\urladdr{http://www.math.univ-paris13.fr/$\sim${}livernet/}

\thanks{The author thanks the Institut Mittag-Leffler (Djursholm, Sweden) 
where this work was carried out during her stay.}

\keywords{$\S$-module, operad, twisted bialgebra, free
  associative algebra, combinatorial Hopf algebra}
\subjclass[2000]{Primary: 18D50 ; Secondary: 16W30, 16A06}
\date{\today}

\begin{abstract} The purpose of this paper is two fold: we study
  the behaviour of the forgetful functor from $\S$-modules to graded
  vector spaces in the context of algebras over an operad and derive
  the construction of combinatorial Hopf algebras. 
  As a byproduct we obtain freeness and cofreeness results for those
  Hopf algebras.

  Let $\Oubli$ denote  the
  forgetful functor from $\S$-modules to graded vector
  spaces. Left modules over an operad $\PP$ are treated as $\PP$-algebras
  in the category of $\S$-modules. We generalize the results obtained by
  Patras and Reutenauer in the associative case to any operad $\PP$: the functor
  $\Oubli$ sends $\PP$-algebras to $\PP$-algebras. If
  $\PP$ is a Hopf operad the functor $\Oubli$ sends Hopf $\PP$-algebras to Hopf $\PP$-algebras. If the
  operad $\PP$ is regular one gets two different structures of Hopf
  $\PP$-algebras in the category of graded vector spaces. We develop
  the notion of unital infinitesimal $\PP$-bialgebras and prove
  freeness and cofreeness results for Hopf algebras built from Hopf
  operads. Finally, we prove that many combinatorial Hopf algebras arise
  from our theory, as it is the case for various Hopf algebras defined on the faces of the permutohedra
  and associahedra.
\end{abstract}

\maketitle

%%%%INTRODUCTION

\section*{Introduction}

An $\S$-module, also named symmetric sequence, is a graded vector
space $(V_n)_{n\geq 0}$ together with a right action of the symmetric
group $S_n$ on $V_n$ for each $n$. The present paper is concerned with
the study, in an operadic point of view, of the forgetful functor $\Oubli$ from
$\S$-modules to graded vector spaces and its applications.

\medskip

The category $\Smod$   of $\S$-modules is a tensor category.
Motivated by the study of homotopy invariants, Barratt introduced the
notion of {\sl twisted Lie algebras} in \cite{Barratt78}, which are Lie
algebras
in the category of $\S$-modules or -- in the operad context --
{\sl left modules} over the operad $\Lie$. More precisely, a twisted
Lie algebra is  an $\S$-module $(L_n)_{n\geq 0}$ together
with a bilinear operation $[,]$ satisfying, for any $a\in L_p,b\in
L_q$ and $c\in L_r$ the relations
\begin{align*}
&[b,a]\cdot\zeta_{p,q}=-[a,b], \\
&[a,[b,c]]+[c,[a,b]]\cdot \zeta_{p+q,r}+[b,[c,a]]\cdot \zeta_{p,q+r}=0,
\end{align*}
where $\zeta_{p,q}$ is the permutation of $S_{p+q}$ given by 
$\zeta_{p,q}(i)=q+i$ if $1\leq i\leq p$ and $\zeta_{p,q}(i)=i-p$ if
$p+1\leq i\leq p+q$. For instance the $\S$-module $(\Lie(n))_{n\geq 0}$
is a twisted Lie algebra for the bracket induced by the operadic
composition
$$\Lie(2)\otimes \Lie(n)\otimes\Lie(m)\rto \Lie(n+m).$$
There exists also a notion of {\sl twisted associative algebras}--
associative algebras in the tensor category $\Smod$-- and
a notion of {\sl twisted associative bialgebras} -- 
associative bialgebras in the category $\Smod$. Stover proved in
\cite{Stover93} that a Cartier-Milnor-Moore
theorem relates the categories of twisted associative bialgebras and
twisted Lie algebras. Following an idea lying in \cite{Stover93}, Patras and Reutenauer proved in
\cite{PatReu04} that two associative bialgebras arise naturally from a twisted associative
bialgebra $(A,m,\Delta)$: 
the {\sl symmetrized bialgebra} $\bar A=(A,\hat m,\bar\Delta)$ and the {\sl cosymmetrized bialgebra} $\hat
A=(A,\bar m,\hat\Delta)$. In \cite{PatSch05} Patras and Schocker
derived from this construction
some known combinatorial Hopf algebras. The first part of this paper is the generalization
of these
constructions to any operad $\PP$. This generalization is performed in two
steps: the first step will focus on the algebra constructions and the second
step on the coalgebra constructions.

\medskip

Given an operad $\PP$, the following notions are the same
\begin{itemize}
\item[] Twisted $\PP$-algebras (see e.g. \cite{LivPat08}),
\item[] Left modules over $\PP$ (see e.g. \cite{Fr04}) and
\item[] $\PP$-algebras in the category of $\S$-modules.
\end{itemize}
Our first question is the following one: given a $\PP$-algebra $M$ in 
$\Smod$ can one endow the graded vector space
$\Oubli(M)$ with a $\PP$-algebra structure? This question 
comes from the observation that $\oplus_n\PP(n)$ is a $\PP$-algebra in 
$\Smod$ but is not {\sl a priori} a $\PP$-algebra in the category of graded
vector spaces. To convice the reader, one can look at the Lie case: if $\PP=\Lie$ then $\oplus_n\Lie(n)$ is a
twisted Lie algebra, and the bracket is not anti-symmetric since the
action of the symmetric group is not trivial. Our first theorem
\ref{T-produit} states that if we apply a symmetrization to the
twisted $\PP$-algebra structure on $M$ then $\Oubli(M)$ is a
$\PP$-algebra. If $\PP=\As$ we recover the definition of the
symmetrized product $\hat m$ of Patras and
Reutenauer. Our second theorem \ref{T-produitr} states that another
product can be defined if the operad $\PP$ is {\sl regular}, that is, $\PP$ is obtained from
a non-symmetric operad tensored by the regular representation of the symmetric group. Then in case $\PP=\As$ we recover the
product $\bar m$ defined by Patras and Reutenauer.

The second step of the construction involves Hopf operads. We define the notion
of Hopf $\PP$-algebras in the category $\Smod$, so that in case
$\PP=\As$ we recover the notion of twisted associative bialgebras. We
develop the analogues of the constructions of Patras and Reutenauer. In theorem
\ref{T-passinghopf} we prove that a  Hopf $\PP$-algebra $(M,\mu_M,\Delta_M)$ in $\Smod$ gives rise
to a  Hopf $\PP$-algebra in $\grVect$, denoted $\bar M=(\Oubli(M),\hat\mu_{\Oubli(M)},\bar\Delta_{\Oubli(M)})$, analogous to the symmetrized bialgebra construction. This is the {\sl
  symmetrized} Hopf $\PP$-algebra associated to $M$.
If $\PP$ is regular then there is another  Hopf
$\PP$-algebra structure  $\hat M=(\Oubli(M),\bar\mu_{\Oubli(M)},\hat\Delta_{\Oubli(M)})$, analogous to the  cosymmetrized bialgebra construction
(see theorem \ref{T-passinghopfr}): this is the {\sl cosymmetrized}
Hopf $\PP$-algebra associated to $M$. Thus the theory developed by Patras and Reutenauer is a consequence of the regularity of the operad $\As$. The example of a
Hopf $\PP$-algebra in $\Smod$ that  we should have in mind throughout
the text is the $\S$-module $\PP$ itself, if $\PP$
is a connected Hopf operad (see section \ref{S-connectedoperad}).

The case $\PP$ regular has another advantage: we define
{\sl unital infinitesimal $\PP$-bialgebras}, which when restricted to  $\PP=\As$, are unital infinitesimal bialgebras as defined by
Loday and Ronco in \cite{LR04}.  The
theorem \ref{T-infinitesimal}  asserts that
$(\Oubli(M),\bar\mu_{\Oubli(M)},\bar\Delta_{\Oubli(M)})$ 
is a unital infinitesimal $\PP$-bialgebra
if $M$ is a Hopf $\PP$-algebra. It is shown in theorem  \ref{T-infinitesimalg}
that the graded vector space $\oplus_n \PP(n)/S_n$ has also a
structure of  unital infinitesimal $\PP$-bialgebra.
These results combined with the theorem of Loday and Ronco (see
\ref{T-LR}) yield the main theorems of our paper, which have some
importance in the study of combinatorial Hopf algebras. First of all if
$\PP=\As$ and $(A,m,\Delta)$ is a twisted associative bialgebra, then the algebra $(A,\bar m,\bar\Delta)$
is a unital infinitesimal bialgebra and consequently is free and
cofree (see theorem \ref{C-fc}). 

%Moreover, theorem \ref{T-isouiP}
%asserts that in case $\PP$ is a Hopf regular multiplicative operad
%that is, endowed with a Hopf operad morphism $\As\rto\PP$, then
%any  unital infinitesimal $\PP$-bialgebra is isomorphic as a unital
%infinitesimal $\PP$-bialgebra to the vector space $T(V)$ where $V$ is the space of primitive elements.

Before going through the applications to combinatorial Hopf algebras,
let us summarize the results in a table.

\bigskip

\begin{tabular}{|p{2cm}|c|c|c|}
\hline
Operad $\PP$ & $\S$-module $M$ & graded vector space $\Oubli(M)$ & Thm\\ \hline
any & $\PP$-alg $(M,\mu)$ & $\PP$-alg
$(\Oubli(M),\hat\mu)$& \ref{T-produit}\\
\hline
regular &  $\PP$-alg $(M,\mu)$ & $\PP$-alg
$(\Oubli(M),\hat\mu)$& \\
& & $\PP$-alg
$(\Oubli(M),\bar\mu)$& \ref{T-produitr}\\
\hline
& & & \\
Hopf & Hopf $\PP$-alg $(M,\mu,\Delta)$ & Hopf $\PP$-alg
$(\Oubli(M),\hat\mu,\bar\Delta)$&\ref{T-passinghopf} \\
\hline
& & & \\
Hopf regular & Hopf $\PP$-alg $(M,\mu,\Delta)$ & Hopf $\PP$-alg
$(\Oubli(M),\hat\mu,\bar\Delta)$& \\
& &  Hopf $\PP$-alg
$(\Oubli(M),\bar\mu,\hat\Delta)$&\ref{T-passinghopfr} \\
& & u.i. $\PP$-bialg $(\Oubli(M),\bar\mu,\bar\Delta)$&\ref{T-infinitesimal} \\
\hline
%& & & \\
%mult Hopf regular &  Hopf $\PP$-alg $(M,\mu,\Delta)$ &
%$(\Oubli(M),\bar\mu,\bar\Delta)\simeq$& \\
%& & $T(\Prim_{\bar\Delta}(\Oubli(M))$&\ref{T-isouiP} \\
%& && \\
%\hline
\end{tabular}

\bigskip

Given a graded vector space $H$, the question is: how does a Hopf
algebra structure arise on $H$? In the last section, we illustrate
with examples that many combinatorial Hopf algebras arise from our
theory. We distinguish two cases. 

If $\PP$ is a Hopf multiplicative operad then $H=\oplus_n \PP(n)$ has
two structures of Hopf algebras: the symmetrized Hopf algebra $\bar H$ is cofree and the
cosymmetrized algebra  $\hat H$ is free. For instance,
we prove that the Malvenuto Reutenauer Hopf algebra arises as the symmetrized Hopf algebra associated to the operad $\As$ and also as the cosymmetrized Hopf algebra associated to the operad $\Zin$ defining Zinbiel algebras. It gives yet another proof for the Hopf algebra of Malvenuto and Reutenauer to be free and cofree independent from its self-duality.
We prove that Hopf algebra structures  on the
faces of the permutohedra given e.g. by Chapoton in \cite{Chapo00a},
Bergeron and Zabrocky in \cite{BerZab05}
and Patras and Schocker in \cite{PatSch05}, arise from the operad
$\CTD$ of commutative tridendriform
algebras defined by Loday in \cite{Lo07}. We deduce also some freeness results from our theory.

In the second case, we assume that there exists a  Hopf multiplicative regular operad $\PP$ such that  
$H=\oplus_n \PP(n)/S_n$. Our theory implies that $H$ is a Hopf $\PP$-algebra hence a Hopf
algebra and is also a unital infinitesimal $\PP$-bialgebra. We prove that
the Hopf algebra of planar trees described by Chapoton in
\cite{Chapo00a}, and the one of planar
binary trees described by Loday and Ronco in \cite{LR98} arise this
way. As a byproduct we obtain freeness results for these Hopf algebras.

\bigskip

The organization of the present paper is as follows. After some
preliminaries in section 1 we define in section 2 the notion
of algebras over an operad $\PP$ in the category of $\S$-modules and in the category of graded vector spaces. 
We explore the structures of $\PP$-algebras on the underlying graded vector space 
of a $\PP$-algebra in $\Smod$. We prove that such a structure always exists and when
the operad is regular one has an additional structure. We compare this result to the one 
obtained by Patras and Reutenauer \cite{PatReu04} 
in the case of twisted associative algebras. In section 3 we study
Hopf operads and the consequences on the underlying graded vector space of a
Hopf $\PP$-algebra in $\Smod$. The section 4 is the study of unital
infinitesimal $\PP$-bialgebras and states the freeness theorems. 
We develop in section 5 the application to combinatorial Hopf
algebras  by means of examples.

\medskip

Throughout the paper, the ground field is denoted by $\kk$ and all vector spaces are $\kk$-vector spaces.

\medskip

\setcounter{tocdepth}{2}
\tableofcontents

%%%%%% SECTION: les S-modules

\section{$\S$-modules and related functors}

%%%%%%%%%% SOUS-SECTION:  Le groupe symetrique

\subsection{The symmetric group}\label{S-symgroup} In this section we develop some material on the symmetric 
group needed in the paper.  The set $\{1,\ldots,n\}$ is written $[n]$. For any set of integers 
$S$, the set $\{p+s,s\in S\}$ is denoted $p+S$. For sets $S\subset
[n]$ and $T\subset[m]$, the set $S\times T$ is the subset $S\cup (T+n)$
of $[n+m]$.

\medskip

Any permutation $\sigma\in S_n$ is written 
$(\sigma_1,\ldots,\sigma_n)$ with $\sigma_i=\sigma(i)$. There is a natural injection
$$\begin{array}{ccc}
S_n\times S_m & \rightarrow &S_{n+m} \\
(\sigma,\tau)&\mapsto &\sigma\times\tau= (\sigma_1,\ldots,\sigma_n,\tau_1+n,\ldots,\tau_m+n).
\end{array}$$

\medskip

The {\sl standardisation} of a sequence of distinct integers $(a_1,\ldots,a_p)$ is the unique 
permutation $\sigma\in S_p$ following the conditions 
$$\sigma(i)<\sigma(j)\Leftrightarrow a_i<a_j,\ \forall i,j.$$
For instance
$$\st(2,13,9,4)=(1,4,3,2).$$

\medskip

Any subset $A=\{a_1<\ldots <a_p\}\subset [n]$ induces
a map
$$\begin{array}{ccc}
S_n & \rightarrow &S_p \\
\sigma &\mapsto &\sigma|_A= \st(\sigma(a_1),\ldots,\sigma(a_p)).
\end{array}$$
For instance
$$(2,6,1,3,5,4)|_{\{1,2,4\}}=\st(2,6,3)=(1,3,2).$$
If $A$ is the empty set then $\sigma|_{\emptyset}=1_0\in S_0$.

\medskip

A {\sl $(p_1,\ldots,p_r)$-shuffle} is a permutation of $S_{p_1+\ldots+p_r}$ 
of type
$$(\tau_1^1,\ldots,\tau_{p_1}^1,\ldots,\tau_1^r,\ldots,\tau_{p_r}^r  )^{-1}$$
with $\tau_1^k<\ldots<\tau_{p_k}^k$ for all $1\leq k\leq r$. The set of all $(p_1,\ldots,p_r)$-shuffles
is denoted by $\Sh_{p_1,\ldots,p_r}$.
For simplicity, a $(p_1,\ldots,p_r)$-shuffle is written as a $r$-uple $(A_1,\ldots,A_r)$ where
$A_1\sqcup\ldots\sqcup A_r$ is an ordered partition of the set $[p_1+\ldots+p_r]$.
Some of the $A_i$'s may be empty.

For instance $(\{2,5\},\{1,3,4\})$ denotes the $(2,3)$-shuffle $(3,1,4,5,2)$. 
Recall that $\Sh_{p_1,\ldots,p_r}$ constitutes a set of right coset representatives
for $S_{p_1}\times\ldots\times S_{p_r}\subset S_{p_1+\ldots+p_r}$, i.e. any 
$\sigma\in S_{p_1+\ldots+p_r}$ has a unique factorization 
$$\sigma=(\sigma_1\times\ldots\times \sigma_r)\alpha,$$
where $\sigma_i\in S_{p_i}$ and where $\alpha$ is a $(p_1,\ldots,p_r)$-shuffle. 
More precisely, if $r=2$
$$\sigma=(\sigma|_{\sigma^{-1}([p])}\times \sigma|_{\sigma^{-1}(p+[q])})(\sigma^{-1}([p]),
\sigma^{-1}(p+[q])).$$

%%%%%%%%%%%%% les S-modules.
\subsection{Graded vector spaces and $\S$-modules} 

%%%%%%%%%%%%%%%%%%%%%definition

\subsubsection{Definition} A {\sl graded vector space} $A$ is a collection $\{A_n\}_{n\geq 0}$ of
$\kk$-vector spaces $A_n$ indexed by the non-negative integers. One can define also $A$ as
$A=\oplus_n A_n$. A map $A\rto B$ of graded vector spaces is a collection of linear morphisms
$A_n\rto B_n$. The category of graded vector spaces is denoted $\grVect$.

An {\sl $\S$-module} $M$ is a graded vector space together with a right 
$S_n$-action $M_n\otimes\kk[S_n]\rto M_n$ for each $n\geq 0$. A map $M\rto N$ of $\S$-modules
is a collection $M_n\rto N_n$ of morphisms of right $S_n$-modules.
The category of $\S$-modules is denoted $\Smod$.

There is a forgetful functor 
$$\Oubli:\Smod\rto\grVect$$ 
which forgets the action of the symmetric group.

%%%%%%%%%%%%%%%%%%% Tensor product

\subsubsection{Tensor product}\label{S-tensorproduct} \hfill\break

The category $\grVect$ 
is a  linear symmetric monoidal category with the 
following {\sl tensor product}:
$$(A\otimes B)_n=\bigoplus_{p+q=n} A_p\otimes B_q.$$
The symmetry isomorphism $\tau:A\otimes B\rto B\otimes A$ is given by
$$\begin{array}{cccc}
\tau:& A_p\otimes B_q &\rto & B_q\otimes A_p \\
& a\otimes b&\mapsto & b\otimes a
\end{array}
$$
The symmetry isomorphism $\tau$ induces a left action of the symmetric group $S_k$ on 
$A^{\otimes k}$, for $A\in\grVect$.

\medskip

The category $\Smod$ is a linear symmetric monoidal category 
with the 
following {\sl tensor product}:
\begin{align*}
(M\otimes N)(n)=&\bigoplus_{p+q=n} (M(p)\otimes N(q))\otimes_{\kk[S_p\times S_q]} \kk[S_n]\\
=&\bigoplus_{p+q=n} (M(p)\otimes M(q))\otimes\kk[\Sh_{p,q}].
\end{align*}
Since a $(p,q)$-shuffle is uniquely determined 
by an ordered partition $I\sqcup J$ of $[p+q]$, 
an element in $(M\otimes N)(p+q)$ can be written $m\otimes n\otimes (I,J)$. 
For the sequel $m\otimes n$ denotes the element $m\otimes n\otimes ([p],p+[q])$
of $(M\otimes N)(p+q)$.
The right action of the symmetric group is given by

\begin{equation}\label{E-sym-partition}
(m\otimes n\otimes (I,J))\cdot \sigma=
m\cdot \sigma|_{\sigma^{-1}(I)}\otimes n\cdot\sigma|_{\sigma^{-1}(J)} 
\otimes (\sigma^{-1}(I),\sigma^{-1}(J))
\end{equation}
The unit for the tensor product is the 
$\S$-module $\un$ given by
$$\un(n)=\begin{cases} \kk, & {\rm if}\ n=0,\\ 0,& \ {\rm otherwise}.
\end{cases}$$
The symmetry isomorphism $\tau:M\otimes N\rto N\otimes M$ is given by
\begin{equation*}
\tau(m\otimes n\otimes (I,J))=n\otimes m\otimes (J,I).
\end{equation*}

For any $\sigma\in S_k$, the symmetry isomorphism induces an isomorphism 
$\tau_\sigma$ of $\S$-modules  from
 $M_1\otimes\ldots\otimes M_k$ to  $M_{\sigma^{-1}(1)}
\otimes\ldots\otimes M_{\sigma^{-1}(k)}$ given by
\begin{multline}\label{F-symmetry}
\tau_\sigma(m_1\otimes\ldots\otimes
m_k\otimes (I_1,\ldots,I_k))=\\
m_{\sigma^{-1}(1)}\otimes\ldots\otimes m_{\sigma^{-1}(k)} \otimes
(I_{\sigma^{-1}(1)},\ldots,I_{\sigma^{-1}(k)}).
\end{multline}

As a consequence $\tau_\sigma$ induces a left $S_k$-action on
$M^{\otimes k}$, for any $\S$-module M.

\medskip

When it is necessary to distinguish the tensor products, we write $\otimes_g$ for the tensor 
product in $\grVect$ and $\otimes_{\S}$ for the one in $\Smod$.

\medskip

The forgetful functor $\Oubli$ does not preserve the tensor product. There are two
natural transformations, $\pi^{\Oubli}$ and $\iota^{\Oubli}$

\begin{align*}
\pi^{\Oubli}_{M,N}, \iota^{\Oubli}_{M,N}: \Oubli(M\otimes_{\S} N)\rto \Oubli(M)\otimes_{g}\Oubli(N)
\end{align*}
defined by, for any $m\in M(p), n\in N(q)$
\begin{align*}
\pi^{\Oubli}_{M,N}(m\otimes n\otimes (I,J))=&\begin{cases}
m\otimes n, & \mbox{\rm if}\ (I,J)=([p],p+[q]), \\
0, & \mbox{\rm otherwise;} \end{cases}\\
\iota^{\Oubli}_{M,N}(m\otimes n\otimes (I,J))= &\ 
m\otimes n.
\end{align*}

When restricted to the full subcategory of vector spaces (the $\S$-modules
concentrated in degree $0$), these two natural transformations restrict to the identity.

%%%%%%%%%%%%%%%%%%%Functors between  $\Smod$ and $\grVect$ 

\subsection{Endofunctors induced by an $\S$-module}

\subsubsection{Endofunctors in $\Smod$} The category of $\S$-modules is endowed with another 
monoidal structure (which is not symmetric): the {\sl plethysm $\circ$}. 
$$(M\circ N)(n):=\bigoplus_{k\geq 0} M(k)\otimes_{S_k} (N^{\otimes k})(n),$$
where $S_k$ acts on the left on
$(N^{\otimes k})$ by formula (\ref{F-symmetry}).
The left and right unit for the plethysm is the $\S$-module $I$ given by
$$I(n)=\begin{cases} \kk, & {\rm if}\ n=1,\\ 0,& \ {\rm otherwise}.
\end{cases}$$

Hence any $\S$-module $M$ defines a functor

$$\begin{array}{cccc}
F_M:& \Smod&\rightarrow &\Smod \\
& N&\mapsto& M\circ N
\end{array}$$
satisfying
$$\begin{cases}
F_I=&\Id \\
F_{M\circ M'}=&F_MF_{M'}   \end{cases}$$

\subsubsection{Endofunctors in $\grVect$}For $M\in\Smod$ and $A\in \grVect$, one can use the same definition
for the plethysm:
$$(M\circ A)(n):=\bigoplus_{k\geq 0} M(k)\otimes_{S_k} (A^{\otimes k})(n).$$
where the tensor product $A^{\otimes k}$ is taken in $\grVect$. Similarly any
$\S$-module $M$ defines a functor

$$\begin{array}{cccc}
F^g_M:& \grVect&\rightarrow &\grVect \\
& A&\mapsto& M\circ A
\end{array}$$
satisfying
$$\begin{cases}
F^g_I=&\Id \\
F^g_{M\circ M'}=&F^g_MF^g_{M'}   \end{cases}$$

\subsubsection{Example}\label{Ex-com} Here is an example that emphasizes the fact that the two functors are 
different even if evaluated at the same underlying vector space. 
Consider the $\S$-module $\Com(n)=\kk$ with the trivial $S_n$-action.
A vector space $V$ is considered either as a graded vector space concentrated in degree
$1$ or as an $\S$-module concentrated in degree 1. This gives
\begin{align*}
F^g_{\Com}(V)=& \oplus_{n\geq 0} \kk\otimes_{S_n} V^{\otimes_g n}=S(V), \\
F_{\Com}(V)=& \oplus_{n\geq 0} \kk\otimes_{S_n} V^{\otimes_{\S} n}=T(V).
\end{align*}

%%%%%%%%%%%%%%%%%%%% Proposition: definition de $\psi$

\subsection{Proposition}\label{P-psi}\it Let $M,N$ be two $\S$-modules. The map
$$\begin{array}{ccc}
M\circ \Oubli(N)& \rightarrow & \Oubli(M\circ N) \\
m\otimes n_1\otimes\ldots\otimes n_k & \mapsto &
\sum_{(T_1,\ldots, T_k)} m\otimes n_1\otimes\ldots\otimes n_k
\otimes(T_1,\ldots,T_k),
\end{array}$$
where $n_i\in N(l_i)$ and $(T_1,\ldots,T_k)$ is an ordered partition of 
$[l_1+\ldots+l_k]$ with $|T_i|=l_i$,
defines a natural transformation
$$\psi_M: F^g_M\Oubli \rightarrow \Oubli F_M$$
functorial in $M\in\Smod$. Furthermore the following diagram commutes
$$\xymatrix{F^g_{M\circ N}\Oubli \ar_{F^g_M\psi_N}[dr]\ar^{\psi_{M\circ N}}[rr]
&& \Oubli F_{M\circ N} \\
& F^g_M\Oubli F_N \ar_{\psi_M F_N}[ur] & }$$ \rm

\medskip

\noindent{\sl Proof.}  Relation (\ref{F-symmetry}) implies that
\begin{multline*}
m\cdot\sigma\otimes n_1\otimes\ldots\otimes n_k\otimes (T_1,\ldots,T_k)=\\
m\otimes n_{\sigma^{-1}(1)}\otimes\ldots\otimes n_{\sigma^{-1}(k)}\otimes(T_{\sigma^{-1}(1)},
\ldots,T_{\sigma^{-1}(k)}).
\end{multline*}
Since the sum is taken over all ordered partitions $(T_1,\ldots,T_k)$, the
image of $m\cdot\sigma\otimes n_1\otimes\ldots\otimes n_k$ is
$$\sum_{(U_1,\ldots, U_k)} m\otimes  n_{\sigma^{-1}(1)}\otimes\ldots\otimes 
n_{\sigma^{-1}(k)}\otimes (U_1,\ldots,U_k)$$
with $|U_i|=l_{\sigma^{-1}(i)}$,
which is the image of $m\otimes n_{\sigma^{-1}(1)}\otimes\ldots\otimes n_{\sigma^{-1}(k)}$.
As a consequence the map is well defined.

The known formula 
\begin{multline}\label{F-shuffles}
\Sh_{p_1^1,\ldots,p_1^{l_1},\ldots,p_k^1,\ldots,p_k^{l_k}}=\\
(\Sh_{p_1^1,\ldots,p_1^{l_1}}\times\ldots\times\Sh_{p_k^1,\ldots,p_k^{l_k}})\Sh_{p_1^1+\ldots+p_1^{l_1},
\ldots,p_k^1+\ldots+p_k^{l_k}} 
\end{multline}
yields the commutativity of the diagram. \hfill $\Box$

%%%%%%%%%%%%%%%%%%%%SOUS-SECTION : LES OPERADES

\section{Algebras over an operad}

In this section we give the definitions of operads and algebras over
an operad and we refer to Fresse \cite{Fr04} for more general theory on operads.
We further state the main results of the section: the underlying
graded vector space of an algebra
over an operad $\PP$ in $\Smod$ is always a $\PP$-algebra and when
$\PP$ is regular there exists a second $\PP$-algebra structure.

%%%%%%%%%%%%%%% def: operad
\subsection{Operads}\label{D-operad} An {\sl operad} is a 
monoid in the category of $\S$-modules
with respect to the plethysm. Namely, an operad is an 
$\S$-module $\PP$ together with a product $\mu_\PP: \PP\circ\PP\rto \PP$ and a unit 
$\eta_\PP: I\rto \PP$ satisfying 
\begin{align*}
\mu_\PP(\Id_\PP\circ\mu_\PP)=&\mu_\PP(\mu_\PP\circ\Id_\PP), \\
\mu_\PP(\Id_\PP\circ\eta_\PP)=&\mu_\PP(\eta_\PP\circ \Id_\PP)=\Id_\PP.
\end{align*}
As a consequence the functors $F_\PP$ and $F^g_\PP$ are monads in the category $\Smod$ and $\grVect$.

The product $\mu_\PP$ is expressed in terms of maps called compositions
$$\begin{array}{ccc}
\PP(n)\otimes\PP(l_1)\otimes\ldots\otimes \PP(l_n)&\rightarrow &\PP(l_1+\ldots+l_n)\\
\mu\otimes\nu_1\otimes\ldots\otimes\nu_n & \mapsto & \mu(\nu_1,\ldots,\nu_n)
\end{array}$$
which are morphisms of right $S_{l_1+\ldots+l_n}$-modules and which factors through the quotient by
the action of the symmetric group $S_n$.

\medskip

The operad $\As$ is the $\S$-module $(\kk[S_n])_{n\geq 0}$.  For $\sigma\in S_n$, 
$\tau_i\in S_{l_i}$
the composition $\mu_{\As}(\sigma;\tau_1,\ldots,\tau_n)$ is the permutation of $S_{l_1+\ldots+l_n}$
obtained by substituting  the block
$\tau_i+l_{\sigma^{-1}(1)}+l_{\sigma^{-1}(2)} \ldots+l_{\sigma^{-1}(\sigma(i)-1)}$ for the integer $\sigma_i$. For instance

$$\mu_{\As}((3,2,1,4);(2,1),(1,3,2),(1),(2,3,1))=(\underbrace{6,5}_{\tau_1+4},\underbrace{2,4,3}_{\tau_2+1},
\underbrace{1}_{\tau_3},\underbrace{8,9,7}_{\tau_4+6}).$$

%%%%%%%%%%%%%%%%%%%%%%Algebras
\subsection{Algebras over an operad}\label{D-algebra} Let $\CC$ denotes either the category of 
$\S$-modules or the category
of graded vector spaces. For any $\S$-module $M$, the functor $F^{\CC}_M$ denotes the functor
$F_M$ or $F^g_M$.

\medskip

Let $\PP$ be an operad. A 
{\sl $\PP$-algebra} or an {\sl algebra over $\PP$} is an algebra over the monad $F^{\CC}_{\PP}$, that is an  object $M$ of $\CC$ 
together with a product
$\mu_M: F^{\CC}_\PP(M)\rto M$ such that the following diagrams commute:

$$\xymatrix{F^{\CC}_{\PP\circ\PP}(M) \ar_{F^ {\CC}_{\mu_\PP}(M)}[d] \ar^>>>>>{F^{\CC}_\PP(\mu_M)}[r]
& F^{\CC}_\PP(M) \ar^{\mu_M}[d] \\
F^{\CC}_\PP(M)\ar^{\mu_M}[r] & M} \qquad \qquad\xymatrix {F^{\CC}_I(M)\ar_{F^{\CC}_{\eta_\PP}(M)}[d]
\ar^{=}[r] & M\\
F^{\CC}_\PP(M) \ar_{\mu_M}[ur] &}$$

\medskip

For $p\in\PP(n)$ and $m_1,\ldots,m_n\in M$ the product $\mu_M(p\otimes m_1\otimes\ldots\otimes m_n)$
is usually written
$$p(m_1,\ldots,m_n)\in M.$$

\medskip

In the category of graded vector spaces one gets the  usual definition of an algebra over an operad.
In the category of $\S$-modules, $\PP$-algebras are also 

-Left modules over $\PP$ in the terminology of Fresse
  \cite{Fr04}, 
  
- If $\PP=\As$ or $\PP=\Lie$, {\sl twisted
  associative or twisted Lie algebras} in the terminology of Barratt
\cite{Barratt78},

- {\sl Twisted $\PP$-algebras} in the terminology of Livernet and Patras
\cite{LivPat08}. 

In the sequel we dedicate the word {\sl twisted} to the only case
$\PP=\As$: a {\sl twisted algebra} is an
algebra over the operad $\As$ in the category $\Smod$.

\medskip

Any free $\PP$-algebra in the category $\CC$ writes $F^{\CC}_{\PP}(M)$ for some $M\in\CC$. As a consequence $\PP$ is
the free $\PP$-algebra in $\Smod$ generated by the $\S$-module $I$.

%%%%%%%%%%%%%%%%%%%%%Theoreme sur les produits

\subsection{Relating $\PP$-algebras in $\Smod$ and in $\grVect$}

%%%%%%%%%%%%%%%%%%le theoreme

\subsubsection{Theorem}\label{T-produit}\it Let $M\in \Smod $ be an algebra over 
an operad $\PP$. The graded vector space
$\Oubli(M)$ is a $\PP$-algebra for the product $\hat\mu_{\Oubli(M)}$
given by the composition
$$\xymatrix{F^g_{\PP}\Oubli(M)\ar^{\psi_{\PP}(M)}[r]& \Oubli F_{\PP}(M)\ar^{\Oubli(\mu_M)}[r]
&\Oubli(M)}.$$ \
That is 
\begin{equation}\label{F-nouveauproduit}
\hat\mu_{\Oubli(M)}(p\otimes m_1\otimes\ldots\otimes m_n)=
\mu_M(p\otimes m_1\otimes\ldots\otimes m_n)\cdot q_{l_1,\ldots,l_n}
\end{equation}
with $p\in\PP(n),m_i\in M(l_i)$
and $q_{l_1,\ldots,l_n}$ is the sum of all $(l_1,\ldots,l_n)$-shuffles. \rm

\medskip

\noindent{\sl Proof.} One has to prove the commutativity of the following two diagrams:

$$\xymatrix{F^g_{\PP\circ\PP}\Oubli(M) \ar_{F^g_{\mu_\PP}\Oubli(M)}[d] 
\ar^>>>>>{F^g_{\PP}\hat\mu_{\Oubli(M)}}[r]
& F^g_\PP\Oubli(M) \ar^{\hat\mu_{\Oubli(M)}}[d] \\
F^g_\PP\Oubli(M)\ar^{\hat\mu_{\Oubli(M)}}[r] & \Oubli(M)} \qquad \qquad\xymatrix {
F^g_I\Oubli(M)\ar_{F^g_{\eta_\PP}\Oubli(M)}[d]
\ar^{=}[r] & \Oubli(M)\\
F^g_\PP\Oubli(M) \ar_{\hat\mu_{\Oubli(M)}}[ur] &}$$
The second diagram is commutative because $\psi_N$ is functorial in $N$, so
\begin{align*}
\psi_I=&\Id, \\
\psi_\PP(F^g_{\eta_\PP}\Oubli)=&(\Oubli F_{\eta_\PP})\psi_I, 
\end{align*}
and because $\mu_M(F_{\eta_\PP}(M))=\Id_M$.

\medskip

Since $M$ is a $\PP$-algebra $\mu_M(F_{\mu_\PP}M)=\mu_M(F_\PP\mu_M).$

The commutativity of the first diagram is a consequence of the computation

$$\begin{array}{ll}
\hat\mu_{\Oubli(M)}F^g_{\mu_\PP}\Oubli(M)& \\
=\Oubli(\mu_M)\psi_\PP(M)(F^g_{\mu_\PP}\Oubli(M))& \mbox{\rm by definition}, \\
=\Oubli(\mu_M)(\Oubli F_{\mu_\PP})\psi_{\PP\circ\PP}(M)& \mbox{\rm by functoriality of\ } \psi, \\
=\Oubli(\mu_M)\Oubli (F_{\PP}\mu_M)\psi_{\PP\circ\PP}(M) & M\ \mbox{\rm is a\ } \PP-{\rm algebra}, \\
=\Oubli(\mu_M)\Oubli (F_{\PP}\mu_M)(\psi_\PP F_\PP)(F^g_\PP\psi_\PP) & \mbox{\rm by proposition \ref{P-psi}},\\
=\Oubli(\mu_M)\psi_\PP(M)(F^g_\PP\Oubli(\mu_M))(F^g_\PP\psi_\PP)\quad  & \psi\ 
\mbox{\rm is a natural transformation},\\
=\hat\mu_{\Oubli(M)}(F^g_\PP\hat\mu_{\Oubli(M)}) & \mbox{\rm by definition}. \hspace{3cm} \Box
\end{array}$$

\medskip

As pointed out in section \ref{D-algebra} any free $\PP$-algebra in $\Smod$ is a $\PP$-algebra and 
satisfies the conditions of the theorem. Hence any free $\PP$-algebra in $\Smod$ gives rise
to a $\PP$-algebra in $\grVect$. In particular 
the graded vector space $\oplus_{n\geq 0} \PP(n)$
is a $\PP$-algebra.

%%%%%%%%%%%%%%%%%%%%%%%%%%%%%%Exemple des algebres commutatives et associatives

\subsubsection{Example}\label{E-comas} We apply formula (\ref{F-nouveauproduit}) for the examples
of the commutative operad and the associative operad.

The commutative operad $\Com$ is the trivial $S_n$-module $\kk$ for all $n$.
Let $e_n$ be a generator of $\Com(n)$. The composition is
$$\mu_{\Com}(e_n\otimes e_{l_1}\otimes\ldots\otimes e_{l_n})=e_{l_1+\ldots+l_n}.$$
The graded vector space $\oplus_n \Com(n)$ is isomorphic to $\kk[X]$. The commutative product on $\kk[X]$ 
induced by the composition $\mu_{\Com}$ is

\begin{equation}\label{E-produitcom}
X^n\hat\cdot X^m={n+m\choose n}X^{n+m},
\end{equation}
since the number of $(n,m)$-shuffles is ${n+m\choose n}$.

\medskip

The associative operad was defined in section \ref{D-operad}. 
The associative product on the space $\oplus_n \kk[S_n]$ induced by
the composition $\mu_{\As}$ is 
\begin{equation}\label{E-produitas}
\sigma\hat*\tau=\sum_{\xi\in\Sh_{p,q}} (\sigma\times \tau)\cdot\xi,
\end{equation}
where $\sigma\in S_p$ and $\tau\in S_q$. This is the product defined by Malvenuto and Reutenauer in 
\cite{MalReu95}.

%%%%%%%%%%%%%%%%%%% Remarque: cas trivial

\subsubsection{Remark}\label{R-trivial} Let $\PP$ be an operad and $M$ be 
a $\PP$-algebra in $\Smod$ such that the action of $S_n$ on $M(n)$ is trivial. 
There is another $\PP$-algebra structure on $\Oubli(M)$ given by
\begin{equation}\label{F-produittrivial}
\mu^{t,g}_{\Oubli(M)}(p\otimes m_1\otimes\ldots\otimes m_n)=\mu_M(p\otimes m_1\otimes\ldots\otimes m_n),
\end{equation}
since the formula (\ref{F-symmetry}) together with the trivial action imply the $S_n$-invariance.

If $\PP=\Com$ then $\kk[X]$ is a commutative algebra for the product
$$X^n\cdot X^m=X^{n+m}.$$

In characteristic $0$ the two commutative products on $\kk[X]$
are isomorphic but it is no more the case in characteristic $p$.

%%%%%%%%%%%%%%%%%%%%%%%Operades regulieres.

\subsection{Regular operads} In this section
we prove that any algebra over a regular operad $\PP$ gives rise to two structures 
of $\PP$-algebra on its underlying graded vector
space. This is the generalization to operads of the result of Patras and Reutenauer in \cite{PatReu04}
in the associative case. Note that this generalization holds only for regular operads.

\medskip

\subsubsection{Definition} The forgetful functor $\Oubli:\Smod\rto \grVect$ has a left adjoint,
the symmetrization functor $\Sym: \grVect\rto \Smod$ which associates
to a graded vector space $(V_n)_n$ the $\S$-module
$(V_n\otimes\kk[S_n])_n$, where the action of the symmetric group is
the right multiplication. An $\S$-module $M$ is {\sl regular} if there
exists a graded vector space $\tilde M$ such that $M=\Sym\tilde M$.
For instance, the $\S$-module $I$ is regular, since $I=\Sym I$. 
Let $\Smodr$ be the subcategory of $\Smod$ of
regular modules (and regular morphisms). A {\sl regular operad} $\PP=\Sym\tilde\PP$ is an operad in the category $\Smodr$, i.e. 
$\mu(\nu_1,\ldots,\nu_k)\in\tilde\PP$ as soon as $\mu,\nu_i\in\tilde\PP$.

Indeed, there is also a plethysm in the category $\grVect$:
$$V\circ^g W=\oplus_k V_k\otimes W^{\otimes k}.$$
Note that
$$\Sym V\circ \Sym W=\Sym(V\circ^g W).$$
A {\sl non-symmetric operad} is a monoid in  the category $\grVect$
with respect to the plethysm $\circ^g$. The operad  $\Sym\tilde\PP$ is
regular if and only if $\tilde\PP$ is a non-symmetric operad.

%%%%%%%%%%%%%proposition: definition de $\psir$
\subsubsection{Proposition}\label{P-psir}\it Let $M=\Sym\tilde M$ be a regular module and $N$ be an 
$\S$-module. The map
$$\begin{array}{ccc}
\tilde M\circ \Oubli(N)& \rightarrow & \Oubli(M\circ N) \\
m\otimes n_1\otimes\ldots\otimes n_k & \mapsto &
m\otimes n_1\otimes\ldots\otimes n_k
\end{array}$$
defines a natural transformation
$$\psi^r_M: F^g_M\Oubli \rightarrow \Oubli F_M$$
functorial in $M\in\Smodr$. Furthermore for $M$ and $N$ in $\Smodr$ the following diagram commutes
$$\xymatrix{F^g_{M\circ N}\Oubli \ar_{F^g_M\psi^r_N}[dr]\ar^{\psi^r_{M\circ N}}[rr]
&& \Oubli F_{M\circ N} \\
& F^g_M\Oubli F_N \ar_{\psi^r_M F_N}[ur] & }$$ \rm

\medskip

\noindent{\sl Proof.}  Since any element in $M$ writes $m\cdot\sigma$ for $m\in\tilde M$, define
$\psi^r_M(m\cdot\sigma\otimes n_1\otimes\ldots\otimes n_k)$ to be
$m\otimes n_{\sigma^{-1}(1)}\otimes\ldots\otimes n_{\sigma^{-1}(k)}$.
It is straightforward to verify the statement with this definition of $\psi^r_M$.
\hfill$\Box$

\medskip

Adapting the proof of theorem \ref{T-produit} by using $\psi^r_{\PP}$
in place of  $\psi_{\PP}$, 
we prove the following theorem:

%%%%%%%%%%%%%%%%%%%%theorem: deuxieme structure de P-algebre

\subsubsection{Theorem}\label{T-produitr}\it Let $M\in \Smod$ be an algebra over a regular operad $\PP$. 
The graded vector space
$\Oubli(M)$ is a $\PP$-algebra for the product $\bar\mu_{\Oubli(M)}$ given by the composition
$$\xymatrix{F^g_{\PP}\Oubli(M)\ar^{\psi^r_{\PP}(M)}[r]& \Oubli F_{\PP}(M)\ar^{\Oubli(\mu_M)}[r]
&\Oubli(M)},$$
that is, for all $p\in\tilde\PP(k)$ and $m_1,\ldots,m_k\in M$,
$$\bar\mu_{\Oubli(M)}(p\otimes m_1\otimes\ldots\otimes m_k)=\mu_M(p\otimes  m_1\otimes\ldots\otimes m_k).$$
Hence $\hat\mu_{\Oubli(M)}$
and $\bar\mu_{\Oubli(M)}$ endow $\Oubli(M)$ with two structures of
$\PP$-algebra.\hfill $\Box$
\rm

\medskip

If $M$ is a trivial $\S$-module  $\bar\mu_{\Oubli(M)}$ coincides 
with $\mu^{t,g}_{\Oubli(M)}$.

Since any free $\PP$-algebra is a $\PP$-algebra, the theorem holds
for any free $\PP$-algebra $F_{\PP}(M)$. In particular 
the graded vector space $\oplus_{n\geq 0} \PP(n)$
is endowed with two structures of $\PP$-algebra.

%%%%%%%%%%%%%%%%%%%Operades multiplicatives

\subsection{Multiplicative operads} An operad $\PP$ is {\sl
  multiplicative} if there exists a morphism of operads
$\As\rto\PP$. Any algebra in $\Smod$ over a multiplicative operad $\PP$
is a twisted algebra and thus its underlying graded vector space is endowed with  two associative
products. It holds in particular for the graded vector space 
$\oplus_n \PP(n)$. 
For instance $\Com$ is a multiplicative operad and we recover the two associative (and commutative)
structures found in example \ref{E-comas} and remark \ref{R-trivial}.

%%%%%%%%%%%%%%%%%%%%%%  SECTION %%%%%%%%%%%%  Left hopf modules and hopf algebras

\section{Hopf algebras over a Hopf operad}

In this section, we generalize the results of Patras and Reutenauer in
\cite{PatReu04} obtained in the associative case to any Hopf operad. We
prove that any Hopf $\PP$-algebra in $\Smod$ yields a
Hopf $\PP$-algebra in $\grVect$  and two Hopf $\PP$-algebras in
$\grVect$ if $\PP$ is regular.

%%%%%%%%%%%%%%%%%%%%% SUBSECTION %%%%%%%%%%%%%%% HOPF OPERADS

\subsection{Hopf operads--the general case}

%In \cite{Moerdijk02}, Moerdijk defined the notion of {\sl Hopf monad} in a tensor category $(\CC,\otimes,U)$. 
%It is a monad $(S,\mu,\eta)$ equipped with maps
%$\tau_{X,Y}: S(X\otimes Y)\rto S(X)\otimes S(Y)$ natural in $X$ and $Y$ and $\theta:S(U)\rto U$
%compatible with the tensor structure of $\CC$ and with the monad structure of $S$.

From now on $\CC$ is either the category $\grVect$ or the category $\Smod$.
The definitions and propositions related to Hopf operads and Hopf algebras over a Hopf operad can be found in \cite{LivPat08}. Here
we recall only what is needed for our purpose.

%%%%%%%%%%%%%%%%% definition: hopf operads

\subsubsection{Definition} \label{D-hopfoperad}Let $\Coalg$ be the category of 
coassociative counital coalgebras, that is vector spaces $V$ endowed with
a coassociative coproduct $\Delta: V\rto V\otimes V$ and a counit $\epsilon:V\rto\kk$.  
A {\sl Hopf operad $\PP$} is an operad 
in $\Coalg$, i.e. $\mu_\PP$ and $\eta_\PP$ are morphisms of coassociative counital coalgebras. 
A Hopf operad amounts to the following data: for each $n$
a coproduct $\delta(n):\PP(n)\rto \PP(n)\otimes\PP(n)$ 
and a counit $\epsilon(n):\PP(n)\rto\kk$ preserving the operadic composition and the action of the symmetric group .
We use Sweedler's notation, that is,
$$\delta(\mu)=\sum_{(1),(2)} \mu_{(1)}\otimes\mu_{(2)}.$$
One has maps in $\Smod$ and in $\grVect$
$$\tau_{M,N}:F_{\PP}(M\otimes N)\rto  F_{\PP}(M)\otimes F_{\PP}(N)$$
and
$$\tau^g_{X,Y}:F^g_{\PP}(X\otimes Y)\rto  F^g_{\PP}(X)\otimes F^g_{\PP}(Y)$$
defined by
\begin{multline}\label{D-tau}
\tau_{M,N}(\mu\otimes m_1\otimes n_1\otimes\cdots\otimes m_k\otimes n_k\otimes (A_1,B_1,\ldots,A_k,B_k))=\\
\sum_{(1),(2)} (\mu_{(1)}\otimes m_1\ldots m_k\otimes \st(A_1,\ldots,A_k))\otimes
(\mu_{(2)}\otimes n_1\ldots n_k\otimes \st(B_1,\ldots,B_k))\\
\otimes (\cup A_i,\cup B_i).
\end{multline}
and
\begin{multline}\label{D-taug}
\tau^g_{X,Y}(\mu\otimes x_1\otimes y_1\otimes\cdots\otimes x_k\otimes y_k)=\\
\sum_{(1),(2)} (\mu_{(1)}\otimes x_1\otimes\ldots\otimes x_k)\otimes 
(\mu_{(2)}\otimes y_1\otimes\ldots\otimes y_k).\end{multline}

As a consequence if $M$ and $N$ are $\PP$-algebras in the category $\CC$, 
then $M\otimes N$ is a $\PP$-algebra for the following product
$$\xymatrix{ F^{\CC}_{\PP}(M\otimes N)\ar^<<<<<{\tau^{\CC}_{M,N}}[r] & F^{\CC}_{\PP}(M)\otimes
 F^{\CC}_{\PP}(N)\ar^<<<<<{\mu_M\otimes\mu_N}[r] & M\otimes N}. $$

%%%%%%%%%%%%%%%%%%%%%%%%%%%% theorem : si P est de Hopf, les P-algebres ont un produit tensoriel

%\subsubsection{Proposition}(\cite{Moerdijk02}) \label{P-hopfhatotimes}\it 
%Let $\PP$ be a Hopf operad. The monads $F_{\PP}$ and $F^g_{\PP}$ are Hopf monads. As a consequence,
%the category of $\PP$-algebras is a tensor category.\rm

\subsubsection{Definition}Let $\PP$ be a Hopf operad. A $\PP$-algebra $M$  is a 
{\sl Hopf $\PP$-algebra} if $M$ is endowed with a coassociative 
coproduct and a counit
$$\Delta_M: M\rto M\otimes M\qquad   \epsilon_M: M\rto\un$$
which are morphisms of $\PP$-algebras. For $\PP=\As$, a  Hopf
$\As$-algebra
is named a {\sl twisted bialgebra}.

%%%%%%%%%%%%%%%%%%%%%theorem: une hopf algebre dans smod donne une hopf algebre dans grvect

\subsubsection{Theorem}\label{T-passinghopf}\it The underlying graded vector space of any Hopf $\PP$-algebra $M$ in $\Smod$ is a 
Hopf $\PP$-algebra in $\grVect$. More precisely, the $\PP$-algebra product on $\Oubli(M)$ is $\hat\mu_{\Oubli(M)}$
and the coproduct
$$\bar\Delta_{\Oubli(M)}:\xymatrix{\Oubli(M)\ar^<<<<<{\Oubli(\Delta_M)}[r]& \Oubli(M\otimes M)
\ar^<<<<<{\pi^{\Oubli}_{M,M}}[r]& \Oubli(M)\otimes\Oubli(M)}$$
is a morphism of $\PP$-algebras. This Hopf $\PP$-algebra is denoted
$\bar M$ and named the {\sl symmetrized Hopf $\PP$-algebra} associated
to $M$.
\rm

\medskip

In particular, if for $m\in M(p)$ one writes
$$\Delta(m)=\sum_{{(1),(2)}\atop {S\sqcup T=[p]}} m_{(1)}^{(S,T)}
\otimes  m_{(2)}^{(S,T)}\otimes (S,T)$$ 
then
\begin{align}\label{F-bardelta}
\bar\Delta(m)=\sum_{(1),(2)}\sum_{i=0}^{p} m_{(1)}^{([i],i+[p-i])}
\otimes  m_{(2)}^{([i],i+[p-i])}
\end{align}

\medskip

\noindent{Proof.} One has to prove that the following diagram is commutative:
$$\xymatrix{ F^g_{\PP}\Oubli(M)\ar_{F^g_{\PP}\bar\Delta_{\Oubli(M)}}[d] 
\ar^{\hat\mu_{\Oubli(M)}}[rr] && \Oubli(M) \ar^{\bar\Delta_{\Oubli(M)}}[dd] \\
F^g_{\PP}(\Oubli(M)\otimes\Oubli(M))\ar_{\tau^g_{\Oubli(M),\Oubli(M)}}[d]&& \\
F^g_{\PP}\Oubli(M)\otimes F^g_{\PP}\Oubli(M)\ar_>>>>>>>>{\hat\mu_{\Oubli(M)}\otimes \hat\mu_{\Oubli(M)}}[rr]&
& \Oubli(M)\otimes\Oubli(M)}$$
Recall that 
\begin{align*}
\bar\Delta_{\Oubli(M)}=&\pi^{\Oubli}_{M,M}\Oubli(\Delta_M), \\
\hat\mu_{\Oubli(M)}=&\Oubli(\mu_M)\psi_{\PP}(M),\\
\Delta_M\mu_M=&(\mu_M\otimes\mu_M)\tau_{M,M}F_{\PP}\Delta_M.
\end{align*}
The functoriality and naturality of $\pi^{\Oubli}$ and $\psi_{\PP}$ imply

\begin{multline*}
\bar\Delta_{\Oubli(M)}\hat\mu_{\Oubli(M)}=\pi^{\Oubli}_{M,M}\Oubli(\Delta_M)\Oubli(\mu_M)\psi_{\PP}(M)\\
=\pi^{\Oubli}_{M,M}\Oubli(\mu_M\otimes\mu_M)\Oubli(\tau_{M,M})\Oubli(F_{\PP}\Delta_M)\psi_{\PP}(M)\\
=(\Oubli(\mu_M)\otimes\Oubli(\mu_M))\pi^{\Oubli}_{F_{\PP}(M),F_{\PP}(M)}
\Oubli(\tau_{M,M})\psi_{\PP}(M\otimes M)F^g_{\PP}\Oubli(\Delta_M).
\end{multline*}
Therefore, the commutativity of the previous diagram follows from the commutativity of
the following diagram
$$\xymatrix{
F^g_{\PP}\Oubli(M\otimes M)\ar_{F^g_{\PP}(\pi^{\Oubli}_{M,M})}[d]\ar^{\psi_{\PP}(M\otimes M)}[rr]&&
\Oubli F_{\PP}(M\otimes M)\ar^{\Oubli\tau_{M,M}}[d] \\
F^g_{\PP}(\Oubli(M)\otimes\Oubli(M))\ar_{\tau^g_{\Oubli(M),\Oubli(M)}}[d]&&
\Oubli(F_{\PP}(M)\otimes F_{\PP}(M))\ar^{\pi^{\Oubli}_{F_{\PP}(M),F_{\PP}(M)}}[d]\\
F^g_{\PP}\Oubli(M)\otimes F^g_{\PP}\Oubli(M)\ar_{\psi_{\PP}(M)\otimes \psi_{\PP}(M)}[rr]&&
\Oubli F^g_{\PP}(M)\otimes\Oubli F^g_{\PP}(M)}$$

Let us compute the composition $R=\pi^{\Oubli}_{F_{\PP}(M),F_{\PP}(M)}\Oubli\tau_{M,M}\psi_{\PP}(M\otimes M)$ 
applied to
$$X=\mu\otimes (m_1\otimes n_1\otimes (A_1,B_1))\otimes\ldots\otimes (m_k\otimes n_k\otimes (A_k,B_k))
\in F^g_{\PP}\Oubli(M\otimes M),$$
where $m_i\in M(l_i)$ and $n_i\in M(r_i)$ and $A_i\sqcup B_i$ is an ordered partition of $[l_i+r_i]$.
$$Y:=\psi_{\PP}(M\otimes M)(X)=\sum_{(T_1,\ldots, T_k)}X\otimes (T_1,\ldots, T_k),$$
where the sum is taken over all ordered partitions of $[l_1+r_1+\ldots+l_k+r_k]$ such that $|T_i|=l_i+r_i$.
For  $U_i=T_i(A_i)$ and $V_i=T_i(B_i)$ one has
$$((A_1,B_1)\times\ldots\times (A_k,B_k))(T_1,\ldots,T_k)=(U_1,V_1,\ldots,U_k,V_k).$$
As a consequence
$$Y=\sum_{(U_1,V_1\ldots,U_k,V_k)}\mu\otimes m_1\otimes n_1\otimes\ldots\otimes
m_k\otimes n_k \otimes (U_1,V_1,\ldots,U_k,V_k),$$
where the sum is taken over all ordered partitions of $[l_1+r_1+\ldots+l_k+r_k]$ such that $|U_i|=l_i, |V_i|=r_i$ 
and $\st(U_i)=A_i, \st(V_i)=B_i$. Set $\bar m=m_1\otimes\ldots\otimes m_k$ and
 $\bar n=n_1\otimes\ldots\otimes n_k$. By the definition of $\tau$ (see (\ref{D-tau})),
\begin{multline*}
Z:=\Oubli\tau_{M,M}(Y)=\sum_{(U_1,V_1\ldots,U_k,V_k)}\\
(\mu_{(1)}\otimes \bar m\otimes\st(U_1,\ldots,U_k))
\otimes (\mu_{(2)}\otimes \bar n\otimes\st(V_1,\ldots,V_k))\otimes (\cup U_i,\cup V_i).
\end{multline*}
Furthermore $\pi^{\Oubli}_{F_{\PP}(M),F_{\PP}(M)}(Z)=0$ except in case $(\cup U_i,\cup V_i)=\Id$.
But 
$$(\cup U_i,\cup V_i)=\Id \Leftrightarrow (A_i,B_i)=\Id, \forall i.$$  
For,
if $(\cup U_i,\cup V_i)=\Id$  then $T_i(A_i)<T_i(B_i)$ and
$A_i<B_i$ which is equivalent to $(A_i,B_i)=\Id$.

As a consequence if $\forall i,\ (A_i,B_i)=\Id$, then
$$R(X)=\sum_{{(U_1,\ldots,U_k)}\atop{(V_1,\ldots,V_k)}}
(\mu_{(1)}\otimes \bar m\otimes (U_1,\ldots,U_k))
\otimes (\mu_{(2)}\otimes \bar n\otimes (V_1,\ldots,V_k)),$$
where the sum is taken over all shuffles $(U_1,\ldots,U_k)$ and $(V_1,\ldots,V_k)$.
If there exists $i$ such that $(A_i,B_i)\not=\Id$, then
$R(X)=0$.

The composition $L=(\psi_{\PP}(M)\otimes \psi_{\PP}(M))\tau^g_{\Oubli(M),\Oubli(M)}F^g_{\PP}(\pi^{\Oubli}_{M,M})$
is easier to compute. First of all, if there exists $i$ such that $(A_i,B_i)\not=\Id$, then
$L(X)=0$. Assume that $(A_i,B_i)=\Id,\forall i$. Then
$$W=\tau^g_{\Oubli(M),\Oubli(M)}F^g_{\PP}(\pi^{\Oubli}_{M,M})(X)=\sum_{(1),(2)}
(\mu_{(1)}\otimes \bar m)\otimes (\mu_{(2)}\otimes \bar n)$$
and $(\psi_{\PP}(M)\otimes \psi_{\PP}(M))(W)=R(X)$. 
Thus $R(X)=L(X),\forall X$ and the diagram is commutative. \hfill $\Box$

%%%%%%%%%%%%%%%%%%%%%%%%%connected operads

\subsubsection{Connected operads and connected Hopf operads}\label{S-connectedoperad} 
An operad is {\sl connected} if $\PP(0)=\kk$
and $\PP(1)=\kk$. Let $1_0$ denote 
the unit of $\kk\in\PP(0)$. If $\PP$ is a connected operad, for any subset $S$ of $[n]$, there 
exists a map
$$\begin{array}{ccc}
\PP(n)&\rto& \PP(|S|) \\
\mu&\mapsto &\mu|_S=\mu(x_1,\ldots,x_n)\end{array}$$
where 
$$\begin{cases}x_i=1_1& \mbox{\rm if}\ i\in S, \\
x_i=1_0& \mbox{\rm if}\ i\not\in S.\end{cases}$$
For $\PP=\As$ one recovers the definition given in section
\ref{S-symgroup} for the symmetric group. 

\medskip

A {\sl connected Hopf operad} is a Hopf operad which is connected and such that 
$\epsilon(0):\kk=\PP(0)\rto\kk$
is the identity isomorphism. As a consequence
$$\epsilon(\mu)=\mu|_{\emptyset}.$$
Recall from \cite{LivPat08} that if $\PP$ is a connected Hopf operad then $\PP$ is a Hopf 
$\PP$-algebra in $\Smod$ for the coproduct given by
\begin{align*}
\Delta(\mu)=&\sum_{{(1),(2)}\atop{S\sqcup T=[n]}} \mu_{(1)}|_S\otimes  \mu_{(2)}|_T\otimes (S,T)\\
=&1\otimes\mu+\mu\otimes 1+\sum_{{S\sqcup T=[n]}\atop{S,T\not=\emptyset}} 
\mu_{(1)}|_S\otimes  \mu_{(2)}|_T\otimes (S,T).
\end{align*}
Indeed the map $\Delta$ is the
unique $\PP$-algebra morphism such that $\Delta(1_1)=1_0\otimes 1_1+1_1\otimes 1_0$.

\medskip

It happens in many examples that $\PP$ is not connected and
$\PP(0)=0$. Nevertheless it is sometimes possible
to define a $\PP$-algebra structure on $(\PP_+\otimes \PP_+)_-$ where 
$$\begin{cases} \PP_+(0)=&\kk \\
\PP_+(n)=&\PP(n),\ n>0 \end{cases}\quad \mbox{\rm and}\
\begin{cases}(\PP_+\otimes \PP_+)_-(0)=&0\\
(\PP_+\otimes \PP_+)_-(n)=&(\PP_+\otimes \PP_+)(n),\ n>0.\end{cases}$$
(see for instance \cite{Lo04Leray}). In that case, $\Delta$ is defined as the unique
$\PP$-algebra map such that $\Delta(1_1)=1_0\otimes 1_1+1_1\otimes 1_0$.
These kind of operads will be treated as connected operads.

\medskip

If $\PP$ is connected, one has two examples of
Hopf $\PP$-algebras in $\grVect$:

\begin{itemize}
\item[$\bullet$] $F^g_{\PP}(V)$ for $V$ in $\grVect$ where the
  product is given by  $F^g_{\mu_\PP}(V)$ and the coproduct is given
  by $F^g_{\Delta}(V)$. 

\item[$\bullet$] For any  $\S$-module $M$ the free $\PP$-algebra
  $F_\PP(M)$ is a Hopf $\PP$-algebra in $\Smod$. The symmetrized Hopf $\PP$-algebra
  $\overline{F_\PP(M)}$ is a Hopf $\PP$-algebra  in $\grVect$ by theorem \ref{T-passinghopf}. 
%the $\PP$-algebra product is given by
%  the composition
%$$\xymatrix{ F^g_\PP\Oubli F_{\PP}(M)\ar^{\psi_{\PP} F_{\PP}(M)}[r]& \Oubli
%  F_{\PP} F_{\PP}(M)
%\ar^{\Oubli F_{\mu_{\PP}}(M)}[r]&\Oubli F_{\PP}(M)}$$ and the
%coproduct is
%given by the composition
%$$\xymatrix{ \Oubli F_{\PP}(M)\ar^{\Oubli F_{\Delta}(M)}[r]& \Oubli
%  (F_{\PP\otimes\PP})(M)}=\xymatrix{\Oubli
%  (F_{\PP}(M)\otimes F_{\PP}(M))
%\ar^{\pi^{\Oubli}_{F_{\PP}(M),F_{\PP}(M)}}[d]\\\Oubli F_{\PP}(M)\otimes
%\Oubli F_{\PP}(M)
%}$$

\end{itemize}

%%%%%%%%%%%%%%%%%%%%%%%%%%%5Cas ou l'operade est reguliere

\subsection{Regular Hopf operads}Let $\PP=\Sym\tilde\PP$ be a regular operad. Assume
$(\PP,\delta)$ is a Hopf operad. The operad $\PP$ is a {\sl regular Hopf operad} if
$\delta(\tilde\PP)\subset \tilde\PP\otimes\tilde\PP$. For instance $\As$ is a regular Hopf operad.
We prove that for any regular Hopf operad a Hopf $\PP$-algebra in the category $\Smod$ 
gives rise to two structures of
Hopf $\PP$-algebra in the category $\grVect$. This is a generalization to regular operads 
of a theorem announced by Stover in \cite{Stover93} and proved by 
Patras and Reutenauer in \cite{PatReu04} in the context of twisted bialgebras. Note 
that the hypothesis regular is needed to obtain two structures of Hopf $\PP$-algebras.

\subsubsection{Theorem}\label{T-passinghopfr}\it Let $\PP$ be a regular Hopf operad. 
Let $(M,\mu_M,\Delta_M)$ be a Hopf $\PP$-algebra in $\Smod$. The product $\bar\mu_{\Oubli(M)}$ together
with the coproduct $\hat\Delta_{\Oubli(M)}$ defined as the composite 
$$\xymatrix{\Oubli(M)\ar^<<<<{\Oubli(\Delta_M)}[r]&\Oubli(M\otimes M)\ar^{\iota^g_{M,M}}[r]&\Oubli(M)\otimes\Oubli(M)}$$
endows $\Oubli(M)$ with a structure of Hopf $\PP$-algebra which is cocommutative if
$(M,\Delta_M)$ is.  This Hopf $\PP$-algebra is denoted
$\hat M$ and named the {\sl cosymmetrized Hopf $\PP$-algebra} associated
to $M$. \rm

\medskip

Note that the coproduct $\hat\Delta_{\Oubli(M)}$ is always defined:
for $m\in M(p)$, if the coproduct in $M$ writes
$$\Delta_M(m)=\sum_{{(1),(2)}\atop {S\sqcup T=[p]}} m_{(1)}^{(S,T)}
\otimes  m_{(2)}^{(S,T)}\otimes (S,T),$$ 
then the induced coproduct in $\Oubli(M)$ writes 
\begin{align}\label{F-hatdelta}
\hat\Delta_{\Oubli(M)}
(m)=\sum_{{(1),(2)}\atop {S\sqcup T=[p]}} m_{(1)}^{(S,T)}
\otimes  m_{(2)}^{(S,T)}.
\end{align}

\medskip

\noindent{Proof.} The proof is similar to the proof of theorem \ref{T-passinghopf}. The commutativity 
of the diagram
$$\xymatrix{ F^g_{\PP}\Oubli(M)\ar_{F^g_{\PP}\hat\Delta_{\Oubli(M)}}[d] 
\ar^{\bar\mu_{\Oubli(M)}}[rr] && \Oubli(M) \ar^{\hat\Delta_{\Oubli(M)}}[dd] \\
F^g_{\PP}(\Oubli(M)\otimes\Oubli(M))\ar_{\tau^g_{\Oubli(M),\Oubli(M)}}[d]&& \\
F^g_{\PP}\Oubli(M)\otimes F^g_{\PP}\Oubli(M)\ar_>>>>>>>>{\bar\mu_{\Oubli(M)}\otimes \bar\mu_{\Oubli(M)}}[rr]&
& \Oubli(M)\otimes\Oubli(M)}$$
is a consequence of the commutativity of
the diagram
$$\xymatrix{
F^g_{\PP}\Oubli(M\otimes M)\ar_{F^g_{\PP}(\iota^{\Oubli}_{M,M})}[d]\ar^{\psi^r_{\PP}(M\otimes M)}[rr]&&
\Oubli F_{\PP}(M\otimes M)\ar^{\Oubli\tau_{M,M}}[d] \\
F^g_{\PP}(\Oubli(M)\otimes\Oubli(M))\ar_{\tau^g_{\Oubli(M),\Oubli(M)}}[d]&&
\Oubli(F_{\PP}(M)\otimes F_{\PP}(M))\ar^{\iota^{\Oubli}_{F_{\PP}(M),F_{\PP}(M)}}[d]\\
F^g_{\PP}\Oubli(M)\otimes F^g_{\PP}\Oubli(M)\ar_{\psi^r_{\PP}(M)\otimes \psi^r_{\PP}(M)}[rr]&&
\Oubli F^g_{\PP}(M)\otimes\Oubli F^g_{\PP}(M)}$$

We first evaluate
$R=\iota^{\Oubli}_{F_{\PP}(M),F_{\PP}(M)}\Oubli\tau_{M,M}\psi^r_{\PP}(M\otimes M)$ 
at 
$$X=\mu\otimes (m_1\otimes n_1\otimes (A_1,B_1))\otimes\ldots\otimes (m_k\otimes n_k\otimes (A_k,B_k))
\in F^g_{\PP}\Oubli(M\otimes M),$$
where $\mu\in\tilde\PP(k), m_i\in M(l_i), n_i\in M(r_i)$ and 
$A_i\sqcup B_i$ is an ordered partition of $[l_i+r_i]$.
Set $\bar m=m_1\otimes\ldots\otimes m_k$ and $\bar n=n_1\otimes\ldots\otimes n_k$ 
and
$$(A_1,B_1)\times\ldots\times (A_k,B_k)=(U_1,V_1,\ldots,U_k,V_k).$$
$$R(X)=\sum_{(1),(2)}(\mu_{(1)}\otimes\bar m\otimes \st(U_1,\ldots,U_k))\otimes
(\mu_{(2)}\otimes\bar n\otimes \st(V_1,\ldots,V_k)),$$
and $\st(U_1,\ldots,U_k)=\Id$ and $\st(V_1,\ldots,V_k)=\Id.$

The composite 
$L:=(\psi^r_{\PP}(M)\otimes \psi^r_{\PP}(M))\tau^g_{\Oubli(M),\Oubli(M)}F^g_{\PP}(\iota^{\Oubli}_{M,M})$
evaluated at $X$ gives
$$L(X)=\sum_{(1),(2)}(\mu_{(1)}\otimes\bar m)\otimes 
(\mu_{(2)}\otimes\bar n).$$
Thus $R(X)=L(X),\forall X$ and the diagram is commutative. 

It is clear that if $\Delta_M$ is cocommutative, so is $\hat\Delta_{\Oubli(M)}$. \hfill $\Box$

\medskip

As a consequence, if $\PP$ is a regular Hopf operad then any Hopf $\PP$-algebra $M$ in $\Smod$ gives rise
to two structures of Hopf $\PP$-algebra in $\grVect$. In particular this result holds for $\oplus_n\PP(n)$
and for the underlying graded vector space of any free $\PP$-algebra.

%%%%%%%%%%%%%%%%Cas des operades multiplicatives

\subsection{Application to multiplicative Hopf operads} In a first step we establish that the 
corresponding Hopf structures in case $\PP=\As$ coincide with the ones discovered by Stover \cite{Stover93} and
proved by Patras and Reutenauer
in \cite{PatReu04}. In a second step we apply the above results to multiplicative Hopf operads.

%%%%%%%%%%%%%%%%%%%%%%%%%%%%Le cas associatif.

\subsubsection{The associative case}\label{S-symco} Recall that the operad $\As$ is a regular Hopf operad. Hence 
the underlying graded vector space of a twisted bialgebra 
is endowed with two structures of Hopf algebra. Let $M$ be a twisted bialgebra with
product $m$ and coproduct $\Delta$.

The Hopf algebra $\bar M=(\Oubli(M),\hat m_{\Oubli(M)},\bar\Delta_{\Oubli(M)})$ is described 
for $a\in M(p), b\in M(q)$ by
\begin{align*}
\hat m_{\Oubli(M)}(a,b)&=m(a,b)\cdot q_{a,b}, \\
\bar\Delta_{\Oubli(M)}(a)=&\sum_{i=0}^p a_{(1)}^{([i],i+[p-i])}\otimes a_{(2)}^{([i],i+[p-i])},
\end{align*}
which is the 
{\sl symmetrized bialgebra} associated to the twisted bialgebra $M$ as
in \cite[proposition 15]{PatReu04}.

The Hopf algebra $\hat M=(\Oubli(M),\bar m,\hat\Delta_{\Oubli(M)})$ 
is described 
for $a\in M(p), b\in M(q)$ by
\begin{align*}
\bar m_{\Oubli(M)}(a,b)&=m(a,b),\\
\hat\Delta_{\Oubli(M)}(a)=&\sum_{S\sqcup T=[p]} a_{(1)}^{(S,T)}\otimes a_{(2)}^{(S,T)},
\end{align*}
which is the {\sl cosymmetrized bialgebra}
associated to the twisted bialgebra $M$ as in 
\cite[definition 8]{PatReu04}.

\medskip

A {\sl multiplicative Hopf operad} is a Hopf operad $\PP$ together with an operad morphism $\As\rto\PP$
which commutes with the Hopf structure. As a consequence any Hopf $\PP$-algebra is a Hopf $\As$-algebra.
The result below is a consequence of the previous sections.

%%%%%%%%%%%%%%%%%%%%%Les deux differentes structure de P-algebre de Hopf quand P est reguliere

\subsubsection{Corollary}\it Let $\PP$ be a multiplicative Hopf operad. The underlying graded vector 
space of any Hopf $\PP$-algebra is endowed with two different structures of Hopf algebra.
\hfill $\Box$
\rm

%%%%%%%%%%%%%%%%%%%%%%%%%%  SECTION:  l'unitaire infinitesimal

\section{Unital infinitesimal $\PP$-bialgebras} 

In this section we give some comparison between $\bar\Delta$ and $\bar\mu$ when the operad is regular,
in view of generalizing the theory of {\sl unital infinitesimal bialgebra} developed by Loday and 
Ronco in \cite{LR06}. This yields the definition of unital infinitesimal $\PP$-bialgebras.   
As a consequence we prove that any Hopf algebra over a multiplicative Hopf operad
is isomorphic to a cofree coassociative algebra. Moreover, if $\PP$ is regular then this isomorphism
respects the $\PP$-algebra structure. We study the associative case in detail.

\bigskip

From now on a connected Hopf operad $\PP$ is given.

%%%%%%%%%%%%%%%%%%%%%%%%%%%sous-section: l'unitaire infinitesimal

\subsection{Unital infinitesimal $\PP$-bialgebras}  In this section,
We prove that
the underlying graded vector space of a Hopf $\PP$-algebra is a unital infinitesimal
$\PP$-bialgebra (theorem \ref{T-infinitesimal}). We prove also in theorem
\ref{T-infinitesimalg} that the same result holds for $F^g_{\PP}(V)$ when 
$V$ is a graded vector space such that $V(0)=0$.

A {\sl connected coalgebra} $M$ in $\Smod$ or $\grVect$ is a coalgebra such that $M(0)=\kk$ and such that 
the counit $\epsilon:\kk=M(0)\rto \kk$ is the identity isomorphism. That is 
for $M\in\Smod$ the coproduct writes
\begin{align*}
\Delta(m)= & 1\otimes m+ m \otimes 1+\sum_{{S\sqcup T=[p]}\atop{S,T\not=\emptyset}}
m_{(1)}^{(S,T)}\otimes m_{(2)}^{(S,T)} \otimes (S,T),
\end{align*}
and for $V\in\grVect$ it writes, $\forall v\in V_r$
\begin{align*}
\Delta(v)= & 1\otimes v+ v \otimes 1+\sum_{{p+q=r}\atop{p,q>0}}
m_{(1),p}\otimes m_{(2),q}.
\end{align*}

%For instance $\PP$ is a connected Hopf $\PP$-algebra in $\Smod$. 
%Any free $\PP$-algebra $F_\PP(M)$ with $M(0)=0$
%is a connected Hopf $\PP$-algebra.

%%%%%%%%%%%%%%%%%%%%%%%%%%%5 def principale

\subsubsection{Definition
}\label{D-uiP}Assume $\PP=\Sym\tilde\PP$ is a regular Hopf operad.  
A {\sl unital infinitesimal $\PP$-bialgebra} $M$ is a $\PP$-algebra in $\grVect$ endowed
with a connected coalgebra structure
$\Delta:M\rto M\otimes M$ satisfying the {\sl infinitesimal relation}:
\begin{multline}\label{E-uiP}
\Delta\mu(m_1,\ldots,m_k)=\sum_{j=1}^k \mu(m_1\otimes 1,
\ldots,m_{j-1}\otimes 1,\Delta(m_j),
1\otimes m_{j+1},\ldots,1\otimes m_k)\\
-\sum_{j=1}^{k-1} \mu(m_1\otimes 1,\ldots,m_{j}\otimes 1,
1\otimes m_{j+1},\ldots,1\otimes m_k),
\end{multline}
for $\mu\in\tilde\PP(k)$.
Note that the operad needs to be regular since the infinitesimal relation is not
$S_k$-equivariant.

\medskip

For instance if $\PP=\As$, a unital infinitesimal $\As$-bialgebra is the definition of 
Loday and Ronco in \cite{LR06} of 
a {\sl unital infinitesimal bialgebra} since the previous relation amounts to
$$\Delta(ab)=\Delta(a)(1\otimes b)+(a\otimes 1)\Delta(b)-a\otimes b.$$

\medskip

Let $M$ be a Hopf $\PP$-algebra in $\Smod$ with $\PP$ regular. 
Theorems \ref{T-passinghopf} and \ref{T-passinghopfr} assert
that the underlying graded vector space of $M$ is endowed with
two structures of Hopf $\PP$-algebras in $\grVect$. One is given by $(\hat\mu,\bar\Delta)$ and the 
other one by $(\bar\mu,\hat\Delta)$. The next theorem explores the relation between
$\bar\mu$ and $\bar\Delta$.

%%%%%%%%%%%%%%%%%%%Theorem: une P-algebre de Hopf fournit une P-alegbre unitaire infinitesimale

\subsubsection{Theorem}\label{T-infinitesimal}\it Let $\PP$ 
be a connected regular Hopf operad and $M$ 
be a connected Hopf $\PP$-algebra in $\Smod$.
The product $\bar\mu:=\bar\mu_{\Oubli(M)}$ and coproduct $\bar\Delta=\bar\Delta_{\Oubli(M)}$ 
endow $\Oubli(M)$ with a structure of unital infinitesimal $\PP$-bialgebra.\rm

\medskip

\noindent{\sl Proof.} Recall that
\begin{align*}
\bar\Delta=&\pi^{\Oubli}_{M,M}\Oubli(\Delta) \\
\bar\mu=&\Oubli(\mu_M)\psi_\PP^r(M).
\end{align*}
Following the proof of theorem \ref{T-passinghopf} one has
$$\bar\Delta\bar\mu=(\Oubli\mu_M\otimes\Oubli\mu_M)\pi^{\Oubli}_{F_{\PP}(M),F_{\PP}(M)}
\Oubli\tau_{M,M}\psi^r_{\PP}(M\otimes M)F^g_{\PP}\Oubli\Delta.$$
Let $X=\mu\otimes m_1\otimes\ldots\otimes m_k\in F^g_{\PP}\Oubli(M)$ with $m_i\in M(h_i)$.
\begin{multline*}
Y:=\psi^r_{\PP}(M\otimes M)F^g_{\PP}\Oubli\Delta(X)=\\
\sum \mu\otimes (m_{1(1)}^{(S_1,T_1)}\otimes m_{1(2)}^{(S_1,T_1)}\otimes(S_1,T_1))\otimes\ldots\otimes
(m_{k(1)}^{(S_k,T_k)}\otimes m_{k(2)}^{(S_k,T_k)}\otimes(S_k,T_k)),
\end{multline*}
where the sum is taken over all ordered partitions $S_i\sqcup T_i$ of $[h_i]$ for all $i$. In order to compute
$\Oubli\tau_{M,M}(Y)$, we write $(S_1,T_1)\times\ldots\times (S_k,T_k)$ as 
$(U_1,V_1,\ldots,U_k,V_k)$ which is an ordered partition of $[h_1+\ldots+h_k]$ and $U_i=|S_i|, V_i=|T_i|$. 
It is obvious that
$\st(U_1,\ldots,U_k)=\Id$ and that the same equality holds for the sequence of $V_i$'s. Furthermore if
$S=\cup U_i$ and  $T=\cup V_i$, then $S=S_1\times\ldots\times S_k$ and $T=T_1\times\ldots\times T_k$.
As a consequence
\begin{multline*}
\Oubli\tau_{M,M}(Y)=\\
\sum_{(S,T)} (\mu_{(1)}\otimes m_{1(1)}^{(S_1,T_1)}\otimes\ldots\otimes
m_{k(1)}^{(S_k,T_k)})\otimes (\mu_{(2)}\otimes m_{1(2)}^{(S_1,T_1)}\otimes\ldots\otimes
m_{k(2)}^{(S_k,T_k)})\otimes(S,T),
\end{multline*}
where the sum is taken over all ordered partitions $(S,T)$ of $[h_1+\ldots+h_k]$ and
where  $S=S_1\times\ldots\times S_k$ and $T=T_1\times\ldots\times T_k$ with $S_i,T_i\subset [h_i]$.
The map $\pi^{\Oubli}_{F_{\PP}(M),F_{\PP}(M)}$ is non zero on an ordered partition $(S,T)$ if and only if 
there exists $r$ such that $S=[r]$. For any such $r$ there exists $j$ such that $S_k=[h_k]$ for $k<j$ and
$S_k=\emptyset$ for $k>j$. Since $M$ is connected 
$m_{(1)}^{\emptyset,[h]}\otimes m_{(2)}^{\emptyset,[h]}=1\otimes m$. 
As a consequence
\begin{multline*}
(\Oubli\mu_M\otimes\Oubli\mu_M)\pi^{\Oubli}_{F_{\PP}(M),F_{\PP}(M)}
\Oubli\tau_{M,M}(Y)=\\
\sum_{j=1}^k\sum_{\alpha=1}^{h_j-1} \mu_{(1)}(m_1,\ldots,m_{j-1},m_{j(1)}^{([\alpha],\alpha+[h_j-\alpha])},
1,\ldots,1)\\ 
\otimes  \mu_{(2)}(1,\ldots,1,m_{j(2)}^{([\alpha],\alpha+[h_j-\alpha])},m_{j+1},
\ldots,m_k)\\
+\sum_{j=0}^k \mu_{(1)}(m_1,\ldots,m_{j},1,\ldots,1)\otimes  \mu_{(2)}(1,\ldots,1,m_{j+1},\ldots,m_k).
\end{multline*}
On the other hand let us compute the right hand side of the equation (\ref{E-uiP}):
\begin{multline*}
\sum_{j=1}^k \mu(m_1\otimes 1,
\ldots,m_{j-1}\otimes 1,\Delta(m_j),
1\otimes m_{j+1},\ldots,1\otimes m_k)\\
-\sum_{j=1}^{k-1} \mu(m_1\otimes 1,\ldots,m_{j}\otimes 1,
1\otimes m_{j+1},\ldots,1\otimes m_k)=
\end{multline*}
\begin{multline*}
\sum_{j=1}^k \mu(m_1\otimes 1,
\ldots,m_{j-1}\otimes 1,\Delta'(m_j),
1\otimes m_{j+1},\ldots,1\otimes m_k)+\\
\sum_{j=1}^k \mu(m_1\otimes 1,
\ldots,m_{j-1}\otimes 1,1\otimes m_j+m_j\otimes 1,
1\otimes m_{j+1},\ldots,1\otimes m_k)\\
-\sum_{j=1}^{k-1} \mu(m_1\otimes 1,\ldots,m_{j}\otimes 1,
1\otimes m_{j+1},\ldots,1\otimes m_k)
\end{multline*}
where 
$$\Delta'(m_j)=\sum_{\alpha=1}^{h_j-1}m_{j(1)}^{([\alpha],\alpha+[h_j-\alpha])}
\otimes m_{j(2)}^{([\alpha],\alpha+[h_j-\alpha])}.$$
Thus the left and right hand sides of the equation (\ref{E-uiP}) are equal 
and the theorem is proved. \hfill $\Box$

\subsubsection{Theorem}\label{T-infinitesimalg}\it Let $\PP=\Sym\tilde\PP$ be a connected regular Hopf operad.
Let $V$ be a graded vector space with $V(0)=0$. The free $\PP$-algebra
in $\grVect$ $F^g_{\PP}(V)$ is a
unital infinitesimal $\PP$-bialgebra. The product is given by the usual product on free $\PP$-algebras
and the coproduct is given for $x=\mu\otimes v_1\otimes\ldots\otimes v_k \in F^g_{\PP}(V)$ with
$\mu\in\tilde\PP(k)$ by
$$\bar\Delta(x)=1\otimes x+x\otimes 1+
\sum_{i=1}^{k-1} (\mu_{(1)}|_{[i]}\otimes v_1\ldots\otimes v_i)\otimes 
(\mu_{(2)}|_{i+[k-i]}\otimes v_{i+1}\ldots\otimes v_k).$$\rm

\medskip

\noindent{\sl Proof.} When $\PP$ is regular $F^g_{\PP}(V)=\oplus_n \tilde\PP(n)\otimes V^{\otimes_g n}$, 
hence it is enough to prove the formula (\ref{E-uiP}) for 
$$\bar\Delta\mu(\nu_1\otimes\bar v_1,\ldots, \nu_k\otimes\bar v_k)$$
with $\mu,\nu_i\in\tilde\PP$ and $\bar v_i\in V^{\otimes l_i}$. 
The computation is straightforward. \hfill$\Box$

%%%%%%%%%%%%%%%%%%%%%%5 On utilise As

\subsection{Rigidity for twisted bialgebras} 
Loday and Ronco proved a  theorem
of rigidity for unital infinitesimal bialgebras.
Recall from \cite{LR06} that the fundamental example of a unital infinitesimal bialgebra is 
given by $T^{fc}(V)=F^g_{\As}(V)$ 
where $V$ is a graded vector space concentrated in degree 1 and where the product 
is given by the concatenation and the coproduct
is given by the deconcatenation. Recall also that for a connected coalgebra $C$, 
with a coproduct $\Delta$ and a counit $\epsilon$, 
the space of primitive elements is defined by
$$\Prim_{\Delta}(C)=\{x\in \Ker\epsilon|\Delta(x)=1\otimes x+x\otimes 1\}.$$
Here is the statement of the theorem:

%%%%%%%%%%%%%%rappel du theoreme de Loday Ronco

\subsubsection{Theorem \cite{LR06}}\label{T-LR} \it Any connected unital infinitesimal Hopf bialgebra
$H$ is isomorphic to $T^{fc}(\Prim(H))$.\rm

\medskip

Let $(A,m,\Delta)$ be a connected twisted bialgebra and $\bar A=(A,\hat m,\bar\Delta)$ the symmetrized bialgebra
and $\hat A=(A,\bar m,\hat\Delta)$  the cosymmetrized bialgebra as in
paragraph \ref{S-symco}. The triple  $(A,\bar m,\bar\Delta)$ is a
unital infinitesimal bialgebra by theorem \ref{T-infinitesimal} and then, by
theorem \ref{T-LR} is isomorphic to
$T^{fc}(\Prim_{\bar\Delta}(A))$. Hence $\hat A$ is a free associative
algebra
and $\bar A$ is a cofree coassociative coalgebra.
Assume furthermore that $\kk$ is of characteristic $0$ and $\Delta$ is
cocommutative. Then $\hat A$ is a cocommutative Hopf algebra, and by
the theorem of Cartier-Milnor-Moore, it is the universal enveloping
algebra of its primitive elements. If each $A_n$ is finite
dimensional, since $\hat A$ is free as an
associative algebra, by lemma 22 in \cite{PatReu04} the space of
primitive elements is a free Lie algebra.

These results are summed up in the following theorem:

%%%%%%%%%%%%%%%%%%%%% Corollaire si P multiplicative on obtient une bigebre unitaire infinitesimale.

\subsubsection{Theorem}\label{C-fc}\it Let   $(A,m,\Delta)$ be a
connected twisted bialgebra. The associated symmetrized bialgebra
$\bar A$ is a cofree
coassociative algebra. The associated  cosymmetrized bialgebra $\hat
A$ is a free associative algebra. 

If $\kk$ is of characteristic 0, if $\Delta$ is
cocommutative and if $A_n$ is finite dimensional for all $n$, there
is an isomorphism of Lie algebras 
$$\Prim_{\hat\Delta}(A)=F_{Lie}^g(\Prim_{\bar\Delta(A)}).$$
This isomorphism is  functorial in A.\hfill$\Box$\rm

\medskip

%%%%%%%%%%%%%%%%%%%%%%%55Exemple pour Com

Using the results of Loday and Ronco we have
improved the results of Patras and Reutenauer. Furthermore, if $\PP$
is a connected multiplicative Hopf operad then it provides connected twisted
bialgebras:
indeed, any Hopf $\PP$-algebra in $\Smod$ is a twisted bialgebra. For
instance $\PP$ and more generally $F_{\PP}(M)$ with $M$ an $\S$-module
such that $M(0)=0$ are connected twisted bialgebras.

\subsubsection{Remark}\label{R-Binfty} If $(A,m,\Delta)$ is a connected twisted
bialgebra then \hfill\break
$(A,\hat m,\bar m,\bar\Delta)$ is a connected
{\sl 2-associative bialgebra} in the terminology of Loday and Ronco in
\cite{LR06},
that is $(A,\hat m,\bar\Delta)$ is a Hopf algebra and $(A,\bar
m,\bar\Delta)$ is a
unital infinitesimal bialgebra. By the structure theorem  in \cite{LR06}, one gets
that $\Prim_{\bar\Delta}(A)$ is a $B_\infty$-algebra and $A$ is the
enveloping 2-as bialgebra of its primitive elements.

\medskip

Assume $\PP$ and $V$ satisfy the conditions of theorem
\ref{T-infinitesimalg}. Assume $\PP$ is multiplicative and
$A=F^g_{\PP}(V)$ is finite dimentional in each degree.
Then $A2=(A^*,^t\Delta,^t\bar\Delta,^t m)$ where $m$ is the associative
product induced by the multiplicative structure of $\PP$ is also a
2-associative bialgebra. If $A2$  is connected then it is  the
enveloping 2-as bialgebra of its primitive elements.

%%%%%%%%%%%%%%%%%%%%%%%%%%%%%% Section: rigidite dans le cas regulier

%\subsection{Rigidity for unital infinitesimal $\PP$-bialgebras} 
%Let $\PP$ be a connected regular Hopf operad which is%multiplicative. Any connected unital infinitesimal $\PP$-bialgebra $(M,\mu_M,\bar\Delta)$ 
%is a unital infinitesimal bialgebra. Applying%the theorem \ref{T-LR} of Loday and Ronco, there is an isomorphism of unital infinitesimal bialgebras
%between $M$ and $T^{fc}(\Prim_{\bar\Delta}M)$. We prove in this section 
%that for any vector space $V$, the space $T^{fc}(V)$ %is endowed with a unique structure of unital infinitesimal $\PP$-bialgebra and that the latter isomorphism is a morphism of
%unital infinitesimal $\PP$-bialgebras.

%%%%%%%%%%%%%%%%%prop: unique structure de P-bialgebre unitaire infinitesimale sur T(V)

%%%%%%%%%%%Le theoreme sur les P-bigebres unitaires infinitesimales.

%\subsubsection{Theorem}\label{T-isouiP}\it Let $\PP$ be a connected regular multiplicative operad. 
%Any connected unital infinitesimal $\PP$-bialgebra $(M,\mu_M,\bar\Delta)$ 
%is isomorphic as a unital infinitesimal $\PP$-bialgebra to 
%$T^{fc}(\Prim_{\bar\Delta}M)$. \rm \hfill$\Box$

\section{Application to combinatorial Hopf algebras}

In this section, we would like to apply our previous results to combinatorial Hopf algebras.
The idea is the following: given a graded vector space $H=\oplus_n H(n)$, 
how does a Hopf algebra structure arise on $H$? We present two cases
coming from the two examples detailed in section \ref{S-connectedoperad}.

\medskip

\noindent{\bf Case 1.} The space $H(n)$ is endowed with a right $S_n$-action. We denote
by $H^{\S}$ the associated $\S$-module. Assume there exists
a connected multiplicative Hopf operad structure on $\PP_H=H^{\S}$. 
From section \ref{S-connectedoperad}, we obtain our first result:
there exists a $\PP_H$-algebra product $\mu$ and a coalgebra coproduct
$\Delta$ such that 
$(H^{\S},\mu,\Delta)$ is a Hopf $\PP_H$-algebra. The graded vector space $H$ has a Hopf 
$\PP_H$-algebra structure which is the symmetrized Hopf $\PP_H$-algebra
$\overline{H^{\S}}$ by theorem \ref{T-passinghopf}.

\smallskip

The second result is a direct consequence of theorem \ref{C-fc}: since
the operad $\PP_H$ is multiplicative, there is a  twisted product 
$m: H^{\S}\otimes H^{\S}\rto H^{\S}$.
As a consequence $(H^{\S},m,\Delta)$ is a twisted bialgebra. The
associated symmetrized Hopf algebra $(H,\hat m,\bar\Delta)$ is cofree and the associated
cosymmetrized Hopf algebra $\bar H=(H,\bar m,\hat\Delta)$ is free. In
case $\Delta$ is cocommutative, under the hypothesis of
theorem \ref{C-fc}, $\bar H$
is the enveloping algebra of the free Lie algebra generated by $\Prim_{\bar\Delta}(H)$.
Furthermore, by remark \ref{R-Binfty} the 2-associative bialgebra
$(H,\hat m,\bar m,\bar\Delta)$ is the 2-associative enveloping
bialgebra of its primitive elements: $\Prim_{\bar\Delta}(H)$ is
endowed with a $B_\infty$-structure.

\smallskip

{\bf Case 1} applies also when $H$ is a free $\PP$-algebra in
$\Smod$ generated by an $\S$-module $M$, with
$\PP$ a multiplicative  Hopf operad and $M(0)=0$.

\medskip

\noindent{\bf Case 2.} Assume $\PP^r_H=\Sym H$ is a connected regular Hopf operad.
The graded vector space $H$ is the free graded $\PP^{r}_H$-algebra
generated by the graded vector space $I$. 
As a consequence, the graded vector space
$(H,\mu,\Delta)$ is a Hopf $\PP^{r}_H$-algebra, where $\mu$ is the
$\PP_H^r$-product and  where
$$\Delta(h)=\sum_{S\sqcup T=[|h|]} h_{(1)}|_S\otimes h_{(2)}|_T$$
comes from the regular Hopf operad $\PP_H^r$. Also $(H,\mu,\bar\Delta)$ with
$$\bar\Delta(h)=\sum_{i=0}^{|h|} h_{(1)}|_{[i]}\otimes h_{(2)}|_{i+[|h|-i]}$$
is a unital infinitesimal
$\PP^{r}_H$-bialgebra by theorem \ref{T-infinitesimalg}. 

\medskip

Again, by remark \ref{R-Binfty}, if $H_n$ is finite dimensional, then
$(H^*,^t\Delta,^t\bar\Delta,^tm)$ in which $m$ is the associative
product, is a $2$-associative bialgebra, which is the 2-associative
enveloping bialgebra of its primitive elements.

\smallskip

{\bf Case 2} applies also when $H$ is a free $\PP$-algebra in
$\grVect$ generated by a graded vector space $V$, with 
$\PP$ a Hopf regular operad and $V(0)=0$.

\medskip

We illustrate by some examples that many combinatorial Hopf algebras arise either from case 1
or from case 2.

%%%%%%%%%%%%%Hopf algebra T(V)

\subsection{The Hopf algebra $T(V)$} Let us apply {\bf Case 1} for
$H=F_{\Com}(V)$ where $V$ is considered as an 
$\S$-module concentrated in degree 1. That is $H=T(V)$. 
As a twisted bialgebra, $T(V)$ is endowed with the concatenation
product and with
the following coproduct 
$$\Delta(v_1\otimes\ldots\otimes v_k)=\sum_{S\sqcup T=[k]} \bar
v|_S\otimes\bar v|_{T}\otimes (S,T)$$
where $\bar v|_S=v_{s_1}\otimes\ldots\otimes v_{s_j}$ for
$S=\{s_1<\ldots<s_j\}$. It is cocommutative.
The symmetrized Hopf algebra structure on $T(V)$ is the shuffle product
together with the deconcatenation, whereas the cosymmetrized Hopf algebra structure on $T(V)$
is the dual structure: the product is the concatenation and the
coproduct is the unshuffle coproduct. In characteristic $0$ it is the
enveloping algebra of the free Lie algebra generated by $V$.

%%%%%%%%%%%%%%%%%%%%%%%%%%%%Quand P=As

\subsection{The Malvenuto-Reutenauer Hopf algebra}\label{S-HMR} This Hopf algebra, denoted $H_{MR}$ has been 
extensively studied
in \cite{MalReu95}, in \cite{NCSF6} under the name of free quasisymmetric functions or in \cite{AguSot05}.
The graded vector space considered is $A=\oplus_n\kk[S_n]$. It is the underlying graded vector
space of the operad $\As$ and {\bf Case 1} applies. Recall that the
operad $\As$ gives rise to a cocommutative twisted bialgebra: 
\begin{align*}
m(\sigma,\tau)=&\sigma\times\tau, \\
\Delta(\sigma)=&\sum_{S\sqcup T=[n]} \sigma|_S\otimes
\sigma|_T\otimes (S,T).\\
\end{align*}
The Hopf algebra
$H_{MR}$ is the symmetrized Hopf algebra
$(A,\hat m,\bar\Delta)$. That is, for $\sigma\in S_n,\tau\in S_m$
\begin{align*}
\hat m(\sigma,\tau)=&\sum_{\xi\in\Sh_{p,q}} (\sigma\times\tau)\cdot\xi, \\
\bar\Delta(\sigma)=&\sum_{i=0}^n \st(\sigma_1,\ldots,\sigma_i)\otimes\st(\sigma_{i+1},\ldots,\sigma_n).
\end{align*}
which is not commutative nor cocommutative.

The cosymmetrized Hopf algebra $\hat A=(A,\bar m,\hat\Delta)$ is given by
\begin{align*}
\bar m(\sigma,\tau)=&\sigma\times\tau, \\
\hat\Delta(\sigma)=&\sum_{S\sqcup T=[n]} \sigma|_S\otimes \sigma|_T.
\end{align*}
The latter Hopf algebra is different from the former one or
its dual since it is a cocommutative Hopf algebra.

\smallskip

From {\bf Case 1} we get that $H_{MR}$ is cofree and that $\hat A$ is
free as an associative algebra: it is generated by the {\sl connected permutations}, the ones which
don't write $\sigma\times\tau$ for $\sigma\in S_n,\tau\in S_m$, $n,m>0$.
In characteristic $0$, $\hat A$ is isomorphic to the enveloping algebra of the free
Lie algebra generated by the connected permutations (compare with theorems 20
and 21 in \cite{PatReu04}). Furthermore, $H_{MR}$ together with $\bar
m$ is a 2-associative bialgebra and it is isomorphic to the 2-associative enveloping
bialgebra
generated by the connected permutations: in \cite{NCSF6} and in \cite{AguSot05} a basis of the space of
primitive elements of $H_{MR}$, indexed by the connected permutations is given.

In paragraph \ref{P-linkMalReu} we prove that $H_{MR}$ is free as an
associative algebra, without using the self-duality of $H_{MR}$.

%%%%%%%%%%%%%%%%%%%%%%%%%%%%%%%%Quand P=Perm

\subsection{Hopf algebra structures on the faces of the permutohedron} Recall that $\Com$ is a Hopf operad.
Let $\overline\Com(n)=\begin{cases}\Com(n) & \textrm{if}\ n>0\\ 0 & \textrm{if}\ n=0\end{cases}$. 
The $\S$-module $\Comp=\As\circ\overline\Com$ has for linear basis the
faces of the $n$-permutohedra. Indeed
\begin{align*}
\As\circ\overline\Com(n)=&\oplus_{k\geq 0}\As(k)\otimes_{S_k}(\overline\Com)^{\otimes k}(n)\\
=&\sum_{(I_1,\ldots,I_k)=[n]}\kk
\end{align*}
where the sum is taken over all {\sl set compositions} (or {\sl ordered set partitions})
$(I_1,\ldots,I_k)$ of $[n]$ such that
$I_j\not=\emptyset,\forall 1\leq j\leq k.$
The action of the symmetric group is given, for $\sigma\in S_n$, by
$$(I_1,\ldots,I_k)\cdot\sigma=(\sigma^{-1}(I_1),\ldots,\sigma^{-1}(I_k)).$$
Chapoton described some Hopf 
algebra structures on the graded vector space $\Oubli(\Comp)$ in
\cite{Chapo00a} and in \cite{Chapo00b},
whereas Patras and Schocker described a \twb\  structure on $\Comp$ in \cite{PatSch05}.
Chapoton described a (differential graded) operad structure on $\Comp$ in \cite{Chapo02b} and Loday 
described a (filtered) one in \cite{Lo07}. 

The aim of this section is to apply our operadic point of view {\bf Case 1} in order
to relate  these structures.

\medskip

A set composition of $[n]$ is written as a word in the alphabet
$\{,\}\cup\{i,1\leq i\leq n\}$. For instance
$(14,2,35)$ is the set composition $(\{1,4\},\{2\},\{3,5\})$ of $[5]$.

%%%%%%%%%%%%%%%%%%%%%%structure d'operade

\subsubsection{Operad structures on the faces of the permutohedron.} Both operads
built by Loday in \cite{Lo07} and Chapoton in \cite{Chapo02b} are quadratic binary operads.
They are generated by the commutative operation represented by the set composition $(12)$
and by the operation represented by the set composition $(1,2)$ in $\Comp(2)$.
Let $w_f=(12)+(1,2)+(2,1)$ and $w_g=(1,2)+(2,1)$.
The composition in the operad $\CTD$ described by Loday is given by the following inductive formula
$$\begin{array}{cc}
(12)(\emptyset,P)=0, &(12)(P,\emptyset)=0, \\
(1,2)(\emptyset,P)=0, & (1,2)(P,\emptyset)=P, 
\end{array}$$
\begin{align*}
(12)(P,Q)=&(P_1\cup Q_1,w_f((P_2,\ldots,P_k),(Q_2,\ldots,Q_l))),\\
(1,2)(P,Q)=&(P_1,w_f((P_2,\ldots,P_k),Q)),
\end{align*}
with $P=(P_1,\ldots,P_k)$ a set composition of $[n]$ and $Q=(Q_1,\ldots,Q_l)$ a set composition
of $[m]$ considered as a set composition of $n+[m]$.  By convention, $w_f(\emptyset,\emptyset)=0$.
The degree of the set composition $P$ is $n-k$. Set compositions
of degree $0$ are in 1-to-1 correspondance with permutations.

The operad $\CTD$ is filtered by the degree of set compositions but not graded. It is not regular and
algebras (in the category of vector spaces) over this operad are named
{\sl commutative tridendriform algebras} by Loday, that is vector
spaces endowed
with a product $\prec$ and a commutative product $\cdot$ satisfying
the relations
\begin{align*}
(x\prec y)\prec z=& x\prec(y\prec z+z\prec y+y\cdot z), \\
(x\cdot y)\prec z=&x\cdot (y\prec z), \\
(x\cdot y)\cdot z=&x\cdot(y\cdot z).
\end{align*}

\medskip

The composition in the operad $\Pi$ described by Chapoton has the same
definition except that $w_f$ is replaced by $w_g$. It is graded by the degree of set compositions. Algebras (in the category of graded vector spaces) 
over $\Pi$ are described in \cite{Chapo02b}.

\medskip

These operads are not connected in the strict sense, 
since the composition with $\emptyset\in\Comp(0)$ is not always defined. 
The equalities involving the emptyset above, are needed for an inductive definition 
and are needed in order to build a coproduct $\Comp\rto\Comp\otimes \Comp$, 
as was explained in the paragraph  \ref{S-connectedoperad}
on connected operads.

%%%%%%%%%%%%%%%%%%Calcul de Delta

\subsubsection{Proposition}\label{P-C-Perm}\it The $\S$-module $\Comp$ is a $\CTD$-Hopf algebra
for the coproduct
$$\Delta(P_1,\ldots,P_k)=\sum_{l=0}^k \st(P_1,\ldots,P_l)\otimes\st(P_{l+1},\ldots,P_k)\otimes 
(\cup_{1\leq j\leq l} P_j,\cup_{l+1\leq h\leq k} P_h).$$
The $\S$-module $\Comp$ is a $\Pi$-Hopf algebra
for the same coproduct.\rm

\medskip

\noindent{\sl Proof.} Let ${\mathcal X}$ denote either the operad $\CTD$ or the operad $\Pi$.
Let $w$ denote either $w_f\in\CTD(2)$ or $w_g\in\Pi(2)$. Note
that for any set composition $P$
\begin{align}\label{F-w}
w(P,\emptyset)=w (\emptyset,P)=P.
\end{align}
One needs first to define
the ${\mathcal X}$-algebra structure on $\Comp\otimes\Comp$. This trick is due to Loday and Ronco:
for $x\in{\mathcal X}(2)$,
\begin{multline}\label{D-coproduct}
x(P_1\otimes P_2,Q_1\otimes Q_2)=\\
\begin{cases} \emptyset\otimes x(P_2,Q_2), & \ \mbox{\rm if}\
P_1=Q_1=\emptyset, \\
x(P_1,Q_1)\otimes w(P_2,Q_2)\otimes (P_1\cup Q_1,P_2\cup Q_2),& \ \mbox{\rm otherwise}.\end{cases}
\end{multline}

The coproduct $\Delta:\Comp\rto\Comp\otimes\Comp$ is the unique ${\mathcal X}$-algebra morphism
mapping $(1)$ to $(1)\otimes\emptyset+\emptyset\otimes(1)$.

\medskip

Let $(I_{l_1},\ldots,I_{l_k})$ be the set composition of $[l_1+\ldots+l_k]$ defined by
$$I_{l_j}=l_1+\ldots+l_{j-1}+[l_j],\forall 1\leq j\leq k.$$ 
Let $n=l_1+\ldots+l_k$. 
We prove the formula for such a set composition by induction on
$k$. If $k=1$, the set composition is just $(n)$. For $n=1$ the formula is proved. If $n>1$ then
$(n)=(12)((1),(n-1))$. By induction one has
\begin{multline*}
\Delta(n)=(12)(\Delta(1),\Delta(n-1))=\\
(12)((1)\otimes\emptyset+\emptyset\otimes (1),
(n-1)\otimes\emptyset+\emptyset\otimes(n-1))= 
(n)\otimes\emptyset+\emptyset\otimes(n),
\end{multline*}
because $(12)(P,Q)=0$ if $P$ or $Q$ is empty and because of relation
(\ref{D-coproduct}) and (\ref{F-w}).

If $k>1$, then $X=(I_{l_1},\ldots,I_{l_k})=(1,2)((I_{l_1}),(I_{l_2},\ldots,I_{l_k})).$
By induction
\begin{multline*}
\Delta(X)=(1,2)(\Delta(I_{l_1}),\Delta(I_{l_2},\ldots,I_{l_k}))=\\
(1,2)(I_{l_1}\otimes\emptyset+\emptyset\otimes I_{l_1},
\sum_{j=1}^{k} (I_{l_2},\ldots,I_{l_j})\otimes(I_{l_{j+1}},\ldots,I_{l_k}))=\\
\emptyset\otimes X+
\sum_{j=1}^{k} (I_{l_1},I_{l_2},\ldots,I_{l_j})\otimes(I_{l_{j+1}},\ldots,I_{l_k})),
\end{multline*} 
because $(1,2)(P,\emptyset)=P$ and $(1,2)(\emptyset,P)=0$ and because of relations 
(\ref{D-coproduct}) and (\ref{F-w}).

For any set composition $P=(P_1,\ldots,P_k)$ of $[n]$, there exists $\sigma\in S_n$ such that
$$(P_1,\ldots,P_k)=(I_{l_1},\ldots,I_{l_k})\cdot\sigma.$$ 
One can choose for $\sigma$ 
the shuffle associated to the set composition P. The coproduct $\Delta$ is a morphism of $\S$-modules.
One gets the conclusion with  formula (\ref{E-sym-partition}). \hfill$\Box$

\medskip

In \cite{LivPat08} we proved that the space of primitive elements with respect to $\Delta$ is a suboperad
of the initial operad. The space of primitive elements is clearly the vector space generated by
the set compositions $(n)$, for $n>0$. Then $\Prim_{\Delta}({\CTD})$ is
the operad $\Com$ (compare with \cite{Lo07}).

%%%%%%%%%%%%%%%%les bigebres twistees associees

\subsubsection{Twisted bialgebras associated to the faces of the permutohedron} The operation
$w_f$ (resp. $w_g$) is associative and commutative. As a consequence, the operads $\CTD$ and $\Pi$ 
are Hopf multiplicative operads and give rise to twisted connected commutative (non cocommutative)
bialgebras $H_f=(\Comp,w_f,\Delta)$ and
$H_g=(\Comp,w_g,\Delta)$. 

\medskip

\hspace{0.6cm} {\bf i) The \twb\  $H_f$.} Patras and Schocker \cite{PatSch05}
defined a \twb\  structure 
on  $\Comp$ denoted $\mathcal T=(\Comp,\star,\delta)$ which is the following. The product $*$ is the concatenation
of set compositions and the coproduct $\delta$ is defined for a set composition $P$ of $[n]$ by
$$\delta(P)=\sum_{S\sqcup T=[n]} \st(P\cap S)\otimes \st(P\cap T)\otimes (S,T),$$
where $(P_1,\ldots,P_k)\cap S=(P_1\cap S,\ldots, P_k\cap S)$ and if $P_i\cap S=\emptyset$ the i-th term
is omitted. For instance $(14,2,35)\cap \{1,3,5\}=(1,35)$. It is clear that $*$ is the dual of $\Delta$
and one can check that $\delta$ is the dual of $w_f$. It gives an operadic interpretation
of their structure: 

\smallskip

\it The dual of the \twb\  defined by Patras and
Schocker is the free commutative tridendriform algebra on one
generator in the category $\Smod$. \rm

\smallskip

Applying {\bf Case 1} one gets that the symmetrized Hopf algebra
$\bar{\mathcal T}$ associated to $\mathcal T$ is cofree, and that the associated cosymmetrized
Hopf algebra $\hat{\mathcal T}$ is free generated by 
{\sl reduced set compositions}: a set composition which is non reduced is the concatenation of two 
non trivial set compositions. For instance $(13,24,6,5)$ is non reduced since it is the concatenation
of $(13,24)$ and $(2,1)$. Moreover if the field $\kk$ is of
characteristic 0, then $\hat{\mathcal T}$ is isomorphic to the enveloping algebra of
the free Lie algebra generated by reduced set compositions.  (Compare with
proposition 10 and corollary 13 in \cite{PatSch05}).
Applying remark \ref{R-Binfty} one gets that $(\mathcal
T,\bar\star,\hat\star,\bar\delta)$ is the 2-associative enveloping
bialgebra on its primitive elements.

\smallskip

The Hopf algebra structure given by Chapoton in \cite{Chapo00a} is
$(\Comp,\bar w_f,\hat\Delta)$ which is the dual of the symmetrized
Hopf algebra $(\mathcal T,\hat\star,\bar\delta)$.  It is also the Hopf algebra $NCQSym$ of
Bergeron et Zabrocky in \cite{BerZab05} and we recover that it is a
free algebra.

\medskip

\hspace{0.6cm} {\bf ii) The \twb\  $H_g$.} This \twb\  gives rise to two Hopf
algebras, which are $(\Comp,\hat w_g,\bar\Delta)$ and $(\Comp,\bar w_g,\hat\Delta)$. One can check
that the latter Hopf algebra is the one described by Chapoton in \cite{Chapo00b}. Again it is a free associative algebra
because $(\Comp,\bar w_g,\bar\Delta)$ is a unital infinitesimal bialgebra. The space of primitive elements
$\Prim_{\bar\Delta}(\Comp)$ is generated by
reduced set compositions. One can check  by an inductive argument, that the Hopf algebra
described by Chapoton is a free associative algebra generated by
reduced set compositions.

%%%%%%%%%%%%%%%%%%%%%%%des set compositions aux permutations

\subsubsection{From set compositions to permutations.}\label{P-linkMalReu} Let
$\Comp_0$ be the sub $\S$-module of $\Comp$ of set compositions of degree
$0$. The vector space $\Comp_0(n)$ is isomorphic to $\kk[S_n]$ but the right $S_n$-action
is given by $\sigma\cdot \tau=\tau^{-1}\sigma$. The $\S$-module
$\Comp_0$ is a sub-operad of $\Pi$. It is the operad $\Zin$,  as noticed by
Chapoton
in \cite{Chapo02b}. In the category of vector spaces, an algebra over $\Zin$ 
is a Zinbiel algebra, that is, a vector
space $Z$ together with a product $\prec$ satisfying the relation
$$(x\prec y)\prec z=x\prec (y\prec z+z\prec y),\ \forall x,y,z\in Z.$$ 
As a consequence there are surjective morphisms of Hopf operads
$$\CTD\rto\Zin, \qquad \Pi\rto\Zin.$$
The operad $\Zin$ is consequently a multiplicative operad and $\Comp_0$
is a commutative twisted
bialgebra. The product and coproduct are given, for $\sigma\in \Comp_0(p)$
and $\tau\in \Comp_0(q)$ by
\begin{align*}
m_{\mathcal Z}(\sigma,\tau)=&\sum_{x\in\Sh_{p,q}} (\sigma\times\tau)x \\
\Delta_{\mathcal Z}(\sigma)=&\sum_{i=0}^p
\sigma|_{[i]}\otimes\sigma|_{i+[p-i]}\cdot (\cup_{1\leq j\leq
  i}\{\sigma(j)\},\cup_{i+1\leq k\leq p}\{\sigma(k)\})
\end{align*}
The morphisms above induce surjective morphisms of twisted bialgebras
$$H_f\rto \Comp_0, \qquad H_g\rto \Comp_0.$$
The cosymmetrized Hopf algebra associated to $\Comp_0$ is clearly
$H_{MR}$, and since it is a cosymmetrized algebra associated to a \twb\
it is free on $\Prim_{\bar\Delta_{\mathcal Z}}(\Comp_0)$. 
But $\bar\Delta_{\mathcal
  Z}(\sigma)=\sum_{\rho\times\tau=\sigma}\rho\otimes\tau.$ 
As a consequence,   $H_{MR}$ is free generated by the connected
permutations and cofree (see section \ref{S-HMR}). We recover the
results obtained in e.g. \cite{PoiReu95}, \cite{NCSF6} and \cite{AguSot05}.

\medskip

Considering the graded linear duals, one has an embedding of cocommutative twisted bialgebras
$$(\Comp_0)^*\hookrightarrow (H_f)^*={\mathcal T}.$$
The symmetrized Hopf algebra associated to $(\Comp_0)^*$ is the dual of
the Malvenuto-Reutenauer Hopf algebra $H_{MR}^*$ (which is isomorphic
to $H_{MR}$).
By functoriality in theorem \ref{C-fc}, we obtain in
characteristic $0$ an embedding of enveloping algebras at the level of
associated cosymmetrized algebras (compare with Theorem 17 in \cite{PatSch05}).

%%%%%%%%%%%%%%%%%%%%%%%%%%%%%%%%%%%SECTON: L'associahedre

\subsection{Hopf algebra structures on the faces of the associahedron}
In his thesis, Chapoton considered various Hopf algebra structures on the faces
of the associahedra, or Stasheff polytopes, filtered in
\cite{Chapo00a}, graded in \cite{Chapo00b}. He considered also
filtered and graded operad structures on these objects in \cite{Chapo02a}.
The filtered operadic structure coincides with the one defined by
Loday and Ronco in \cite{LR04}, under the name of {\sl tridendriform
  operad}. In this section,
we apply {\bf Case 2} to obtain Hopf algebra structures on the
faces of the associahedra.

\medskip

\subsubsection{Planar trees} The set of planar trees with $n+1$
leaves is denoted by $T_n$. The set $T_n=\cup_{k=0}^{n-1} T_{n,k}$
is graded by $k$ where $n-k$ is the number of internal vertices. For instance
$T_{n,0}$ is the set of planar binary trees. The Stasheff polytope of
dimension $n-1$ has its faces of dimension $0\leq k\leq n-1$ indexed by $T_{n,k}$.
The aim of this section is to provide the vector space $\oplus \kk[T_n]$
with Hopf structures.

\medskip

Given some planar trees $t_1,\ldots,t_k$ 
the planar tree $\vee(t_1,\ldots,t_k)$ is the one obtained by joining
the roots  of
the trees $t_1,\ldots,t_k$ to an extra root, from left to right.
If $t_i$ has degree $l_i$ then $\vee(t_1,\ldots,t_k)$ has degree
$l_1+\ldots+l_k+k-1$.

One can label the $n$ sectors delimited by a tree $t$ in $T_n$ from left
to right as in the following example:

$$t=\epsfxsize=2in\epsfbox{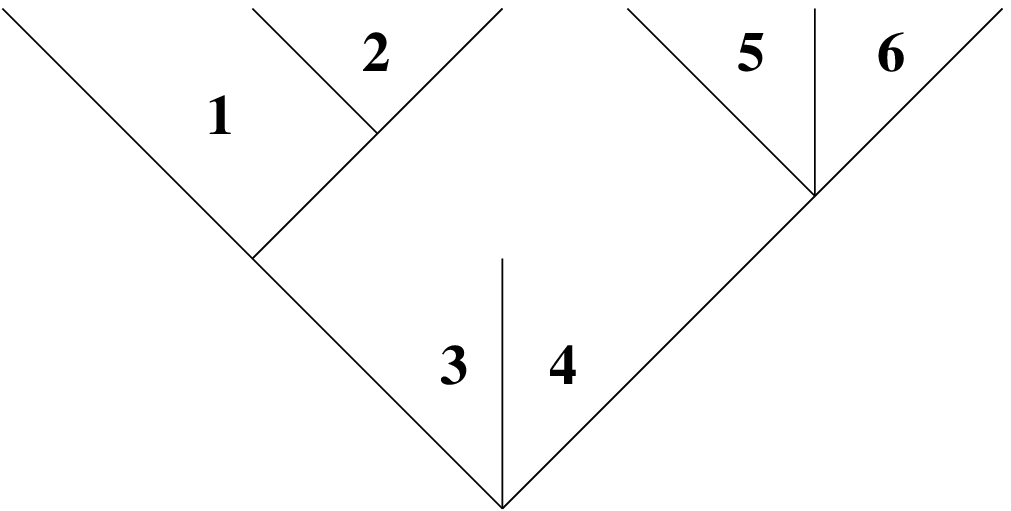}=\vee(\epsffile{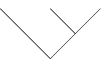},\epsffile{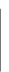},\epsffile{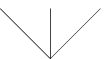}).$$

%%%%%%%%%%%%%%%%%%%%%%%%%%%%%%%Les algebres tridendriformes

\subsubsection{The operad of tridendriform algebras \cite{LR04}, \cite{Chapo02a}} The operad $\TD$
is a regular operad whose underlying $\S$-module is
$\TD(n)=\Sym\kk[T_n]$.
It is a quadratic binary operad generated by $3$ operations
\begin{align*}
\prec:=& \epsffile{1_2.eps}, \\
\succ:=& \epsffile{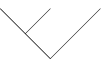}, \\
\cdot:=& \epsffile{12.eps},
\end{align*}
satisfying the relations
$$\begin{array}{l}
\begin{cases}
(x\prec y)\prec z=x\prec(y*z)=& \epsffile{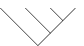}+\epsffile{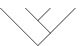}+\epsffile{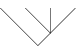}, \\
(x\succ y)\prec z=x\succ(y\prec z)=&\epsffile{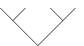}, \\
(x*y)\succ z=x\succ(y\succ z)=& \epsffile{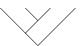}+\epsffile{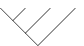}+\epsffile{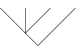},
\end{cases}\\
\begin{cases}
(x\succ y)\cdot z=x\succ(y\cdot z)=&\epsffile{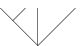},\\
(x\prec z)\cdot z=x\cdot (y\succ z)=&\epsffile{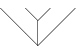},\\
(x\cdot y)\prec z=x\cdot(y\prec z)=&\epsffile{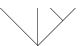},
\end{cases} \\
\begin{cases} (x\cdot y)\cdot z=x\cdot(y\cdot z)=&\epsffile{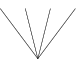},\end{cases}
\end{array}$$
where
\begin{align}\label{E-*} 
x*y=&x\cdot y+x\prec y+x\succ y
\end{align}
is associative.

One can also give an inductive formula for the composition in $\TD$ by
the following
\begin{align*}
 \epsffile{1_2.eps}(|,y)=y, &\quad \epsffile{1_2.eps}(x,|)=0, \\
 \epsffile{12.eps}(|,y)=0, &\quad \epsffile{12.eps}(x,|)=0, \\
  \epsffile{2_1.eps}(|,y)=0, &\quad \epsffile{2_1.eps}(x,|)=x, \\
 \epsffile{1_2.eps}(\vee(x_1,\ldots,x_k),y)=&\vee(x_1,\ldots,x_{k-1},x_k*y), \\
 \epsffile{12.eps}(\vee(x_1,\ldots,x_k),\vee(y_1,\ldots,y_l))=&\vee(x_1,\ldots,x_{k-1},x_k*y_1,y_2,\ldots,y_l), \\
 \epsffile{2_1.eps}(x,\vee(y_1,\ldots,y_l))=&\vee(x*y_1,y_2,\ldots,y_l). \\
\end{align*}
In  \cite{Chapo02a}, Chapoton describes the composition $x\circ_i y$
for trees $x$ and $y$: 
$$x\circ_i y=\sum_{(f_l,f_r)} x\circ_i^{(f_l,f_r)} y,$$
where $y$ is inserted in the sector $i$ of $x$ following the maps
$(f_l,f_r)$: the left (right) edge of the sector $i$ of $x$ is a
set of edges and vertices ordered from bottom to top and denoted 
by $x^i_l$, ($x^i_r$). The left (right) most edge of $y$ has several
vertices: the ordered set of these vertices is denoted by $y_l$ ($y_r$).
The map $f_l$ ($f_r$) is an increasing map from $y_l$ ($y_r$) to
$x^i_l$ ($x^i_r$). For instance
$$\epsfxsize=1.2in\epsfbox{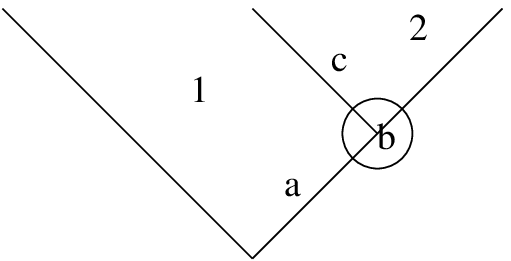}\circ_1 \epsfxsize=1.2in\epsfbox{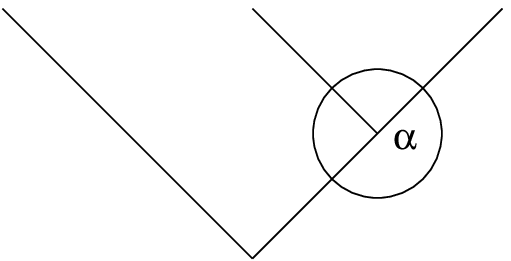}=
  \underbrace{\epsffile{1_2_3.eps}}_{f_1}+\underbrace{\epsffile{1_3_2.eps}}_{f_2}+
\underbrace{\epsffile{1_23.eps}}_{f_3},$$
where 

$$\begin{array}{c}
x^1_l=y^l=\emptyset, \\
x^1_r=\{a<b<c\}, \ y_r=\{\alpha\}, \\
f_1=(\Id,\alpha\mapsto a), \ f_2=(\Id,\alpha\mapsto c)\ \  \hbox{\rm and}
\ \ f_3=(\Id,\alpha\mapsto b).
\end{array}$$

\medskip

Moreover one can define a Hopf structure on this operad following the same pattern than proposition \ref{P-C-Perm}:
for $x\in\TD(2)$,

$$x(t_1\otimes t_2,s_1\otimes s_2)=
\begin{cases} |\otimes x(t_2,s_2), & \ \mbox{\rm if}\
t_1=s_1=|, \\
x(t_1,s_1)\otimes t_2*s_2,& \ \mbox{\rm otherwise}.\end{cases}$$
By induction on the degree of a tree $t$ one can prove that

%%%%%%%%%%%%%%%%%%%%%%%%%%%Calcul de Delta

\subsubsection{Proposition}\it The $\S$-module
$(\kk[T_n]\otimes\kk[S_n])_n$ is a tridendriform Hopf algebra for the
coproduct
\begin{multline*}
\Delta(\vee(t_1,\ldots,t_k))=|\otimes \vee(t_1,\ldots,t_k) +\\
\sum_{S_i\sqcup T_i=[l_i]}
\vee(t_{1(1)}^{S_1},\ldots,t_{k(1)}^{S_k})\otimes
t_{1(2)}^{T_1}*\ldots *t_{k(2)}^{T_k}\otimes 
(S_1\bar\times\ldots\bar\times S_k,T_1\bar\times\ldots\bar\times
T_k)
\end{multline*}
where $|t_i|=l_i$, 
$$\Delta(t_i)=\sum_{S_i\sqcup T_i=[l_i]} t_{i(1)}^{S_i}\otimes
t_{i(2)}^{T_i}\otimes (S_i,T_i)$$
and for any $U_i\subset[l_i]$,
$$U_1\bar\times\ldots\bar\times U_k=\{l_1+\ldots+l_{i}+i, 1\leq i\leq
k-1\}\cup_{1\leq i\leq k} U_i+l_1+\ldots+l_{i-1}+i-1.\ \Box$$
\rm

%%%%%%%%%%%%%%%%%%%%%%%%%%%%%%%Les structures de Hopf

\subsubsection{Hopf structures}\label{P-HopfTriDend} Let us apply {\bf Case 2} to the
connected Hopf operad $\TD$
which is multiplicative with the product $*$ introduced in equation (\ref{E-*}).

The graded vector space
$\Tree=\oplus_n\kk[T_n]$ is the free graded tridendriform algebra over one generator,
hence it is a Hopf tridendriform algebra in $\grVect$. Let $\mu$ be
the tridendriform product and $\Delta$ be the coproduct.
By theorem \ref{T-infinitesimalg}, $(\Tree,\mu,\bar\Delta)$
is a unital tridendriform bialgebra. The description of $\Delta$
gives the description of $\bar\Delta$: if
$\bar\Delta(t_k)=\sum t_{k(a)}\otimes t_{k(b)}$, then
$$\bar\Delta(\vee(t_1,\ldots,t_k))=|\otimes \vee(t_1,\ldots,t_k) +\sum
\vee(t_1,\ldots,t_{k-1},t_{k(a)})\otimes
t_{k(b)}.$$
As a consequence a basis of $\Prim_{\bar\Delta}(\Tree)$ is given by
the planar trees of type $\vee(t_1,\ldots,t_{k-1},|)$.

\medskip

\it The Hopf algebra $(\Tree,*,\Delta)$ is the free associative algebra
spanned by the set of trees of the form $\vee(t_1,\ldots,t_{k-1},|)$. \rm

\medskip

Recall that the product $*$ is defined by induction: for $x=\vee(x_1,\ldots,x_k)$ and
$y= \vee(y_1,\ldots,y_l)$
\begin{multline*}
x*y=
\vee(x_1,\ldots,x_{k-1},x_k*y)+\vee(x*y_1,y_2,\ldots,y_l)+\\
\vee(x_1,\ldots,x_{k-1},x_k*y_1,y_2,\ldots,y_l).
\end{multline*}
and
$$\Delta(t)=|\otimes t+\sum 
\vee(t_{1(1)},\ldots,t_{k(1)})\otimes
t_{1(2)}*\ldots *t_{k(2)}.$$
The Hopf structure defined by Chapoton in \cite{Chapo00a} is essentially the same: the product is the same 
and the coproduct is $\tau\Delta$, where $\tau$ is the symmetry isomorphism. Hence it is
a free associative algebra spanned by the set of trees of the form $\vee(|,t_2,\ldots,t_{k})$.

\medskip

The graded linear
dual of $\Tree$ is a $2$-associative bialgebra: it is free as
an associative algebra for the product $^t\bar\Delta$, cofree as a
coalgebra for $^t*$ and it is the enveloping $2$-as bialgebra of its
primitive elements (with respect to $^t*$).
The product of two trees $^t\bar\Delta(t,s)$ is the tree obtained by gluing the
tree $s$ on the right most leave of $t$. it is usually denoted
$t\backslash s$.
For instance
$$\epsffile{1_2.eps}\backslash \epsffile{12.eps}=\epsffile{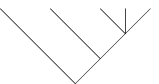}.$$

\subsubsection{Some operad morphisms}
There are morphisms of Hopf operads

\begin{align}\label{DiagramOperad}
\xymatrix{\TD\ar_{\pi_{\TD}}[d]\ar^{\psi}[r]&
  \CTD\ar^{\pi_{\CTD}}[d]\\
\Dend\ar^{\psi_0}[r]& \Zin}
\end{align} 
The vertical maps are projection onto cells of degree $0$. The map
$\pi_{\CTD}$
has been explained in paragraph \ref{P-linkMalReu}. The map
$\pi_{\TD}$ is the projection onto the dendriform operad, which is a
regular operad generated by planar binary trees (see e.g. \cite{LR98}).
The morphism $\psi$
sends $\prec$ to the set composition $(1,2)$, $\succ$ to the set composition $(2,1)$ and
$\cdot$ to the set composition $(12)$. Indeed we can describe $\psi$ at the level of
trees. There is a map $\phi$ from set compositions to trees, described by
induction as follows. Let $P=(P_1,\ldots,P_k)$ be a set composition of $[n]$. 
If $P_1=\{l_1<\ldots<l_j\}$ then it splits $[n]$ into $j+1$
intervals $I_s$ possibly empty: for $0\leq s\leq j,\ I_s=]l_{s},l_{s+1}[$ with $l_0=0$ and $l_{j+1}=n+1$.
The map $\phi$ is defined by 
$$\begin{cases}
\phi(\emptyset)=& |, \\
\phi(P)=& \vee(\phi(P\cap I_0),\ldots,\phi(P\cap I_j)).\end{cases}$$ 

For instance if $P=(34,1,56,2)$ then $I_0=\{1,2\}, I_1=\emptyset$ and $I_2=\{5,6\}$ and
$$\phi(34,1,56,2)=\vee(\phi(1,2),\phi(\emptyset),\phi(12))=\epsfxsize=1.2in\epsfbox{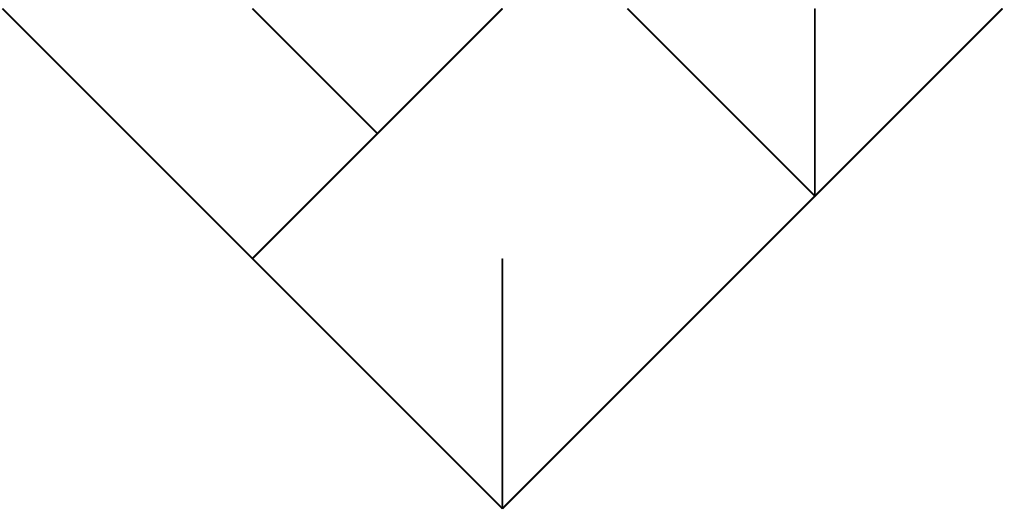}.$$
Note that the function $\theta$ from the permutohedra to the 
associahedra defined by Tonks in \cite{Tonks97} (see also
\cite{Chapo00a} and \cite{PatSch05}) 
satisfies 
$$\theta(P_1,\ldots,P_k)=\phi(P_k,P_{k-1},\ldots,P_2,P_1).$$

The morphism $\psi$ is the transpose of $\phi$.
It is an operad morphism, whereas the transpose of $\theta$ is not an operad morphism.

The morphism $\psi_0$ is the transpose of $\phi_0$ which is an operad morphism.
Loday and Ronco defined in \cite{LR98} a 
function from $\kk[Y_n]$ (the vector space generated
by planar binary trees with $n$ vertices) to $\kk[S_n]$ in order to embed $\kk[Y_n]$ as
a Hopf subalgebra of the (graded linear dual of the) Malvenuto-Reutenauer Hopf algebra. It is also a transpose of a set
morphism $S_n\rto Y_n$. If $\alpha:S_n\rto S_n$ is the involution defined by
$$\alpha(\sigma_1,\ldots,\sigma_n)=(\sigma_n,\ldots,\sigma_1)^{-1},$$
then the set morphism defined by Loday and Ronco is 
$\phi_0\alpha$. 

\subsubsection{Consequences on Hopf algebra morphisms}\hfill\break

%\medskip

The operad $\TD$ is a regular operad. A {\sl tridendriform 2-bialgebra} is a 4-uple
$(H,\mu,\Delta,\delta)$ where $(H,\mu,\Delta)$ is a Hopf tridendriform bialgebra
in $\grVect$ and $(H,\mu,\delta)$ is a unital infinitesimal tridendriform bialgebra.
The diagram (\ref{DiagramOperad}) is a diagram of tridendriform 2-bialgebras:
\begin{itemize}
\item[$\bullet$] The operad $\TD$ induces the tridendriform 2-bialgebra
$(\Tree:=\oplus_n \kk[T_n],\mu_T,\Delta_T,\bar\Delta_T)$ 
explained in paragraph \ref{P-HopfTriDend}. 
\item[$\bullet$] Using the surjective operad morphism $\Pi_{\TD}:\TD\rto \Dend$ one gets the 
tridendriform 2-bialgebra structure
on the vector space spanned by planar binary trees denoted by
$$(PBT:=\oplus_n \kk[Y_n],\mu_{Y},\Delta_{Y},\bar\Delta_{Y}),$$
where $Y_n$ is the set of planar binary trees with $n+1$ leaves as in \cite{LR98}. 
\item[$\bullet$] The underlying $\S$-modules of the operads $\CTD$ and $\Zin$ are
tridendriform algebras in $\Smod$ then by theorems \ref{T-passinghopf} and \ref{T-infinitesimal} one gets
tridendriform 2-bialgebras on the underlying graded vector spaces. These structures are denoted respectively
$(\Comp,\mu_C,\Delta_C,\bar\Delta_C)$ and $(\oplus_n\kk[S_n],\mu_S,\Delta_S,\bar\Delta_S)$.
\end{itemize}
As a consequence we obtain a diagram of Hopf algebras (and unital infinitesimal bialgebras as well):

\begin{align*}
\xymatrix{(\oplus_n \kk[T_n],*_T,\Delta_T)\ar_{\pi_{\TD}}[d]\ar^>>>>>>>>{\psi}[rr]&&
 (\Comp,*_C,\Delta_C)=NCQSym \ar^{\pi_{\CTD}}[d]\\
(\oplus_n\kk[Y_n],*_Y,\Delta_Y)\ar^{\psi_0}[rr]&&  (\oplus_n\kk[S_n],*_S,\Delta_S)=H_{MR}}
\end{align*} 
where the horizontal arrows are injective morphisms of Hopf algebras and
vertical arrows are surjective morphisms of Hopf algebras.
Note that the Hopf algebra structure on the planar binary tree $(\oplus_n\kk[Y_n],*_Y,\tau\Delta_Y)$
where $\tau$ is the symmetry isomorphism is the one described 
by Loday and Ronco in \cite{LR98}. Note also that the graded linear dual of this diagram is a diagram of
$2$-associative bialgebras. It extends the results obtained by Palacios and Ronco in \cite{PalRon06}.

\subsubsection{Conclusion} For the last decade, many results of
freeness and cofreeness of combinatorial Hopf algebras have appeared in the
litterature (see the references cited throughout the paper and recently
\cite{AguSot06}, \cite{Foissy07}, \cite{NovThi06}). The present paper illustrates that
these freeness results are a consequence of an operadic structure on
the Hopf algebra $H$ itself or its symmetrization $\Sym H$. Namely, either
the Hopf algebra $H$ is an $\S$-module and one can find an operad
structure on $H$ in order
to apply {\bf Case 1}; or the Hopf algebra is not an $\S$-module and
one can find an operad structure on $\Sym H$ in order to apply {\bf Case 2}.

\bibliographystyle{plain} 
\bibliography{bibliomai08.bib}

\bigskip

\end{document}